\documentclass{article} 
\usepackage{amsmath,amssymb}
\usepackage{mathrsfs}
\usepackage{amsthm}
\usepackage[dvipdfmx]{graphicx,xcolor}
\usepackage{indentfirst}		

\makeatletter
    
    \@addtoreset{equation}{section}
\makeatother

\newcommand\thmcite{}
	\newtheoremstyle{mystyle}
  		{}
 		{}
		{\rm}
		{0pt}
		{\bfseries}
		{\vspace{0.3\baselineskip}}
 		{\newline}
		{\thmname{#1}~\thmnumber{#2}\thmnote{~(#3)}\thmcite}
 	\theoremstyle{mystyle}
		\newtheorem{dfn}{Definition}[section]
		\newtheorem{thm}{Theorem}[section]
		\newtheorem{lem}{Lemma}[section]
		\newtheorem{cor}{Corollary}[section]
		\newtheorem{prp}{Proposition}[section]
		\newtheorem{exm}{Example}[section]
		\newtheorem{prf}{Proof :}
		
		\def\@endtheorem{%
  			\hfill □
 		 	\endtrivlist}

\usepackage[top=30truemm,bottom=30truemm,left=30truemm,right=30truemm]{geometry}
\usepackage{mathrsfs}

\DeclareMathOperator*{\capp}{\cap}
\DeclareMathOperator*{\cupp}{\cup}
\DeclareMathOperator*{\sqcupp}{\sqcup}

\def\ca{\displaystyle\capp}
\def\cu{\displaystyle\cupp}
\def\sqcu{\displaystyle\sqcupp}

\def\wi{Wiener }				
\def\wii{Wiener-It\^o }			
		
\def\wifs{Wiener functionals }		\def\wifS{Wiener functionals}
\def\ni{Nisio }				\def\nI{Nisio}
\def\sk{Skorokhod }			\def\sK{Skorokhod}
\def\og{Ogawa }				\def\oG{Ogawa}
\def\ue{Uemura }				\def\uE{Uemura}
				
\def\sci{Schwarz inequality }		\def\scI{Schwarz inequality}
				
\def\sos{Sobolev space }			
		\def\sosS{Sobolev spaces}	
\def\ba{Banach }				
\def\bo{Borel }				
			\def\leB{Lebesgue}

\def\nc{noncausal }	\def\nC{noncausal}	\def\Nc{Noncausal }

\def\sfci{SFC }				
\def\sfcis{SFCs }				\def\sfciS{SFCs}

\def\sfco{SFC-O }			\def\sfcO{SFC-O}
\def\sfcos{SFC-Os }			\def\sfcoS{SFC-Os}

\def\sfcs{SFC-S }				\def\sfcS{SFC-S}
\def\sfcss{SFC-Ss }			\def\sfcsS{SFC-Ss}

\def\sfcou{SFC-$\rm{O}_{u}$ }			\def\sfcoU{SFC-$\rm{O}_{u}$}
\def\sfcop{SFC-$\rm{O}_{\vp}$ }		\def\sfcoP{SFC-$\rm{O}_{\vp}$}

\def\sft{SFT }	\def\sfT{SFT}

\def\invt{invertibility }		

\def\i{\item}			%

\def\r#1{\ref{#1}}

\def\tr#1{Theorem \r{#1}}
\def\ssr#1{Subsection \r{#1}}

\def\ip#1#2#3{\langle#1,#2\rangle_{#3}}

\def\blf#1#2#3#4{{}_{#1}\langle#2,#3\rangle_{#4}}

\def\EX#1{E\Bigl(\,#1\,\Bigl)}

\def\aex#1{E|\,#1\,|}

\def\AEx#1{E\Bigl |\,#1\,\Bigl |}

\def\ic#1{\mathsf{1}_{#1}}

\def\tsum{\textstyle\sum\limits}			%

\def\tsup{\textstyle\sup\limits}			%

\def\tinf{\textstyle\inf\limits}			%

\def\tlim{\textstyle\lim\limits}			%

\def\lims#1{\underset{#1}{\varlimsup}}	

\def\limi#1{\underset{#1}{\varliminf}}

\def\tprod{\textstyle\prod\limits}			%

\def\mab#1{\mathbb{#1}}			%

\def\mac#1{\mathcal{#1}}			%

\def\C{\mab{C}}
\def\R{\mab{R}}
\def\Q{\mab{Q}}
\def\Z{\mab{Z}}
\def\N{\mab{N}}

\def\F{\mac{F}}

\def\L{\mac{L}}
\def\M{\mac{M}}

\def\P{\mac{P}}

\def\s#1#2#3{(#1_#2)_{#2\in#3}}

\def\t{\times}			%

\def\inc{\hookrightarrow}			%

\def\om{\omega}			%
\def\Om{\Omega}			%

\def\ot{\otimes}			%

\def\a{\forall}			%

\def\e{\exists}			%

\def\vp{\varphi}			%

\def\BM{Brownian motion }			%

\def\remi#1{\noindent\textbf{Remark #1}\ }

\def\rem{\noindent\textbf{Remark }\ }

\def\sf{$\sigma$-field }		%

\def\du{d_{\rm u}}

\def\sfc#1#2{(#1,#2)}

\def\dub{\,d_{\rm u}B}	%

\def\dpy{d_{\vp}Y}

\def\de{\delta}

\def\d{\displaystyle}

\def\q{\quad}		%

\def\tx{\text}		%

\def\m{\medskip}

\def\bt{ }
\def\et{\@endtheorem}

\def\mpn{\medskip\par\noindent}

\def\ni{\noindent}

\def\inum{\hfill
\vspace{-\abovedisplayskip}
\vspace{-\baselineskip}}

\title{Identification of \nc finite variation processes from the stochastic Fourier coefficients
    } 

\author{%
  Kiyoiki Hoshino\footnote{Osaka Prefecture University, Japan.
    E-mail: su301032@edu.osakafu-u.ac.jp
    }
  }	

\date{}
\begin{document}
\maketitle

{\bf Abstract}
Let $\s{B}{t}{[0,\infty)}$ be a Brownian motion on a probability space $(\Om,\F,P)$.
Our concern is whether and how a \nc type stochastic differential $dX_t=a(t,\om)\,dB_t+b(t,\om)\,dt$ is identified from its stochastic Fourier coefficients (\sfcis for short) $\sfc{e_n}{dX}:=\int_{0}^L\overline{e_n(t)}\,dX_t$ with respect to a CONS $(e_n)_{n\in\mab{N}}$ of $L^2([0,L];\mab{C})$.
This problem has been studied by S.\og and H.\ue (\og (2013)\cite{Og13}, (2014)\cite{Og14}; \og and \ue (2014)\cite{OU14a}, \cite{OU14b}, (2015)\cite{OU15}).
In this paper we give several results on the problem 
for each of stochastic differentials of \og type and \sk type when $[0,L]$ is an finite or infinite interval.
Specifically,
we first give a condition for a random function to be identified from the SFCs and apply it to obtain affirmative answers to the question with several concrete reconstruction formulas of the random functions.
This paper restates the result given in \cite{H19a}
by a metamathematical notion of constructiveness we introduce here.

\section{Introduction}

\indent Let $\s{B}{t}{[0,1]}$ be a one-dimensional \BM defined on a probability space $(\Om,\F,P)$.
It has been discussed by S.\og and H.\ue (\cite{Og13}-\cite{OU17a}) whether and how a 
random function on $[0,1]$, namely, a  complex-valued jointly measurable map $a(t,\om)$ on $[0,1]\t\Om$ is identified from a set of its stochastic Fourier coefficients (\sfcis for short) defined by $\hat{a}_n=\sfc{e_n}{a\,dB}:=\int_0^1\overline{e_n(t)}a(t,\om)\,dB_t$, where $\overline{e_n(t)}$ denotes the complex conjugate of $e_n(t)$, with respect to a CONS $(e_n)_{n\in\mab{N}}$ of $L^2([0,1];\mab{C})$. 
Also, it is discussed in \cite{OU14b} as an extension of the question above whether and how random functions $a(t)$ and $b(t)$ are identified from a set of \sfcis defined by $\sfc{e_n}{dX}:=\int_0^1\overline{e_n(t)}a(t,\om)\,dB_t+\int_0^1\overline{e_n(t)}b(t,\om)\,dt$ of the stochastic differential $dX_t=a(t)\,dB_t+b(t)\,dt$.
Note that the symbol $\int \,dB$ stands for some sort of stochastic integral with respect to $\s{B}{t}{[0,1]}$.
The SFC is called of \sk type if it is defined by the \sk integral (\cite{Sk}, see also Definition \ref{dfn:92}) and of \og type if it is defined by the \og integral (\cite{Og79}, see also Definition \r{def:0001}).

The notion of the \sfci and \sft (stochastic Fourier transformation) were introduced by \og in a series of studies \cite{Og86}-\cite{Og91} related to a Fredholm type stochastic integral equation (SIE) for random fields.
In these articles, the \invt of the \sft is used to show the existence and uniqueness of solutions for SIEs. 
In the previous studies \cite{Og13}-\cite{OU17a}, 
affirmative answers to these questions are given.
In \cite{Og13} and \cite{Og14}, the random functions are causal, i.e. adapted to a filtration for $\s{B}{t}{[0,1]}$.
In \cite{OU14a} and \cite{OU14b}, the random functions are \nC (anticipative), i.e. not necessarily causal and square integrable \wifs
and the \sfcis are of \sk type.
In \cite{OU15} and \cite{OU17a}, the random functions are \nc and non-negative absolutely continuous and the \sfcis are of \og type, and the CONS is the exponential system $(\exp{(2\pi\sqrt{-1}nt)})_{n\in\mathbb{Z}}$
(\tr{thm: 0001}).
On the other hand, T.Kazumi and the author showed  in \cite{KH17b} the identifications of \nc square integrable \wifs and stochastic differentials (as extensions of those) by using the \wii decomposition, where the \sfcis are of \sk or \og type
(\tr{thm: 0003}, \r{thm: 12} and \r{thm: 21}). Concrete statements of Theorems \r{thm: 0001}-\r{thm: 21}, which are related to this note, are collected in \ssr{sub: 1.1}.

Also, in \cite{OU17a} \og and \ue proposed two meanings of identifying a random function from its \sfciS:
the wide sense and the strong sense.
The former simply indicates unique determinability of a random function from its \sfciS.
This is equivalent to \invt of the \sft under proper conditions as in \cite{Og14}.
The latter is derived from an application viewpoint.
Consider \sfcis as given data, you cannot use information of the underlying Brownian motion to estimate quantity such as volatility in finance.
It can be said that the identification problem in the strong sense asks, roughly speaking, whether it is possible to identify a diffusion coefficient or a drift term of a differential $constructively$ and $without$ $the$ $Brownian$ $motion$ from \sfciS. 
The study on the latter is required in 
the study of the volatility estimation problem proposed by P. Malliavin et al. (\cite{MM},\cite{MT}) and conducted by \og and \uE(\cite{Og14},\cite{OU15},\cite{OU17a}).
In this note, in relation to the latter 
we introduce the notion of constructive identification in
an assigned first-order language.
The reason why we introduce the metamathematical notion of constructiveness in Appendix A
is to give a framework to explain and evaluate derivations (or derivation formulas) or any other maps (or formulas of maps)
by an organized notion based on the definite criteria.
In addition, as another type of identification, for purely mathematical interest 
we introduce B-dependent (resp. B-independent) identification,
which can be called identification ''in need of'' (resp. ''in no need of'')
the condition that the underlying Brownian motion is $\s{B}{t}{[0,\infty)}$,
from the purely mathematical interest.
Here, B-dependent (resp. B-independent) is short for dependent on Brownian motion (resp. independent of Brownian motion).
Note that 
for each derivation map of the random function, there could be various different formulas which represent the map,
while these notions of B-dependent and B-independent identifications depend only on the derivation map in itself and does not depend on each derivation formula which represents the derivation map. 

Main aim of this paper is to give the following identifications,
while providing reconstruction formulas, of random functions from \sfciS.\m

\noindent $\cdot$ Identification from \sfcis of \og type (\sfcoS)
\begin{enumerate}
\i[$\cdot$]\tr{thm:6.20}, \tr{pr:4.11466}, Corollary \r{thm: 1.5} (extension of \tr{thm: 0001}) : Constructive identification for any \nc finite variation process.
\i[$\cdot$]\tr{thm:6.220}, \tr{pr:4.11440}, Corollary \r{thm: 1.555} (extensions of \tr{thm:6.20}, \tr{pr:4.11466}, Corollary \r{thm: 1.5}, respectively) : Constructive identification for a stochastic differential whose diffusion coefficient is any \nc finite variation process.
\i[$\cdot$]\tr{thm: a5} (additional result, extension of \tr{thm: 0003}) : Constructive identification of a stochastic differential whose diffusion coefficient is an S-type It\^o process or a more general \wi functional.
\end{enumerate}

\noindent $\cdot$ Identification from \sfcis of \sk type (\sfcsS)
\begin{enumerate}
\i[$\cdot$]\tr{thm:6.2000}, \tr{pr:4.11400}, Corollary \r{thm: 15} : Constructive identification for a stochastic differential whose diffusion coefficient is a locally absolutely continuous \wi functional.
\end{enumerate}
Here in each theorem listed above, the CONS which defines \sfcis is taken generally
and strictly speaking there includes metamathematical assertions since the notion of constructiveness is metamathematical.

The organization of this paper is as follows.
In Section 2, first, we describe fundamental notions concerning the issues, in particular, the \nc stochastic integrals: \og and \sk integrals. 
In Section \r{ssec:3.1001},
we give propositions regarding the \og integrability, which give sufficient conditions for \sfcos to be defined and  representations of \sfcoS, 
which are used to show the results listed above.
In Section 4, first, we give the precise definitions of \sfcis of stochastic differentials.
Next, we introduce definitions of constructive identification in assigned first-order languages
and B-independent identification
of a random function from the \sfciS.
In Section 5, we give several lemmas.
In Sections 6 and 7, 
we first give a necessary and sufficient condition, under proper conditions, for a random function or some quantity to be identified from \sfcis and apply it to give the main results listed above.
In Appendix A, we introduce the notion of constructiveness in an assigned first-order language,
which is used to define constructive identification in Section 4.
In Appendix B, we discuss some measurability and continuity of stochastic processes
to obtain the main results.

\subsection{Preceding studies}
\label{sub: 1.1}

In this subsection, we review preceding results on the identification problems.

\subsubsection{Identification from \sfcoS}

First, we overview results about the identification from \sfcoS.
The following is an extension of the result originally presented in \cite{OU15}.

\begin{thm}[\og\!\!, \ue\tx{\cite{OU17a}}(2018)]
Assume a CONS of $L^2[0,1]$ is given by the exponential system $(T_n)_{n\in\mab{Z}}$ defined by $T_n(t)=\exp{(2\pi\sqrt{-1}nt)}$
and a real-valued random function $a(t)$ on $[0,1]$ satisfies the following:
\begin{enumerate}
\i $P(\,a(t)$ is non-negative and absolutely continuous in $t\,)=1$.
\i $a'(t)\in L^2([0,1]\times\Omega)$.
\i $\int_0^1a(t)\,dt\in L^2(\Omega)$.
\end{enumerate}
Then, we have
\begin{align}
P\Biggl(\,\lims{h\searrow 0}\frac{\mathscr{F}_{\tau, T}(a)'(t+h)-\mathscr{F}_{\tau, T}(a)'(t)}{\sqrt{2h \log\log(\frac{1}{h})}}=a(t)\q
\a t\in[0,1]\,\Biggl)=1,
\label{al:0.00100}
\end{align}
where $\mathscr{F}_{\tau, T}(a)(t,\om)=\hat{a}_0(\om)+\tsum_{n\neq 0}\frac{1}{2\pi\sqrt{-1}n}\hat{a}_n(\om)T_n(t)$ and $(\hat{a}_n)_{n\in\mab{Z}}$ is given by u-integrals for $L^2([0,1];\mab{R})$
(see Definition \r{def:0001}).
\label{thm: 0001}
\end{thm}

In Theorem \ref{thm: 0001} Ogawa and Uemura also showed that $a(t)$ can be identified even from $(\hat{a}_n)_{n\in\Lambda}$, where $\Lambda$ is a cofinite subset of $\mathbb{Z}$.\m

On the other hand, from Theorems 5.1,5.2 and 5.3 in \cite{KH17a}, we can see the following fact about the integrability of the Ogawa integral:
Let $\L_1^{r,2}$ denote the \sos of square integrable \wifs with differentiability index $r\in[1,\infty)$ (see Definition \r{def:00111} for the precise definition).
Given $e:[0,1]\to\mab{R}$ and a CONS $(\vp_m)_{\in\mab{N}}$ of $L^2([0,1];\mab{R})$ and a \sk integral process
\begin{align}
a(t):=\int_0^tf(s)\,\delta B_s\,\,,t\in[0,1],
\label{al:000}
\end{align}
where $\int \,\de B$ stands for the \sk integral and $f\in\L_1^{2,2}$.
We introduce the following conditions (C.1),(C.2) and (C.3) on $e(t)$, $(\vp_m)_{\in\mab{N}}$ and $f(t)$.

\noindent (C.1)
\begin{enumerate}
\i $e(t)$ is of bounded variation.
\i $\tsup_{M\in\mab{N}}\Bigl|\tsum_{m=1}^{M}\vp_m\widetilde{\vp}_m\,\Bigl|_{L^2[0,1]}<\infty$,
where $\widetilde{\vp}_m(t)=\int_0^t\vp_m(s)\,ds$.
\end{enumerate}
\noindent (C.2)
\begin{enumerate}
\i $e(t)$ is regulated,\\
namely, $e(t)$ has finite left and right limits on $(0,1]$ and $[0,1)$, respectively.
\label{en:001}
\i $\tsup_{M\in\mab{N}}\Bigl|\tsum_{m=1}^{M}\vp_m\widetilde{\vp}_m\,\Bigl|_{L^1[0,1]}<\infty$.
\i there exists $g\in\mathcal{L}_1^{1,2}$ such that $f(t)=\int_0^tg(s)\,\de B_s$.
\end{enumerate}
\noindent (C.3)
\begin{enumerate}
\i $e(t)$ is of bounded variation.
\i there exists $g\in\mathcal{L}_1^{1,2}$ such that $f(t)=\int_0^tg(s)\,\de B_s$.
\end{enumerate}
Now, the following assertions (A),(B) and (C) about the \og integrability hold:
\begin{enumerate}
\i[(A)] \, If (C.1) or (C.2) is satisfied, then $eX$ is $\varphi$-integrable.
\i[(B)] \, If (C.3) is satisfied, then $eX$ is u-integrable for $L^2([0,1];\mab{R})$.
\i[(C)] \, If the assumption of (A) or (B) holds, the integral of $eX$ converges in $L^2(\Om)$ and it is given  by 
$$
\int_0^1e(t)a(t)\,\de B_t+\frac{1}{2}\int_0^1e(t)f(t)\,dt
+\int_0^1e(t)\Bigl(\,\int_0^tD_tf(s)\,\de B_s \,\Bigl)\,dt.
$$
\end{enumerate}

\rem The condition \r{en:001} of (C.2) is weaker than ''$e(t)$ is piecewise continuous'' which we assumed in Theorem 5.3 in \cite{KH17a}.
But above assertion is justified by the argument in the proof of Lemma 4.3 in \cite{KH17a} because we can check the space $BV([0,1])$ of functions  of bounded variation is dense in the space ${\rm Reg}([0,1])$ of regulated functions  equipped with the uniform norm (Theorem 7.6.1 in \cite{Dd}).

From (A),(B),(C) and Theorem \r{thm: 21} mentioned below, we have the following similar to Theorem 4.3 in \cite{KH17b}.

\begin{thm}
Let $(e_n)_{n\in\mathbb{N}}$ and $(\vp_m)_{\in\mab{N}}$ be CONSs of $L^2([0,1];\mab{R})$ and $a(t)$ a process defined by \eqref{al:000} and $b\in\mathcal{L}_1^{0,2}$.
Assume that $e_n(t)$, $(\vp_m)_{\in\mab{N}}$ and $f(t)$ satisfy any of the
conditions (C.1),(C.2) and (C.3) for each $n\in\mab{N}$.
Let
$$
d_{\dagger}Y_t:=\begin{cases}
	a(t)\,d_{\rm{u}}B_t+b(t)\,dt& \text{, if (C.3) holds for each $e_n(t)$}\\
	a(t)\,d_{\varphi}B_t+b(t)\,dt& \text{, otherwise,}
	\end{cases}
$$
where $\int\, d_{\rm{u}}B$ and $\int\, d_{\varphi}B$ stand for the u-integral for $L^2([0,1];\mab{R})$ and $\vp$-integral, respectively (see Definition \r{def:0001}).
Then, $a(t)$ and $b(t)$ are identified from the system $((e_n,d_{\dagger}Y))_{n\in\mab{N}}$ of SFC-$\rm{O}_{\dagger}$'s.
\label{thm: 0003}
\end{thm}

\subsubsection{Identification from \sfcsS}

Next, we overview results about the identification from \sfcsS.
The following are two results stated in \cite{KH17b}, which show the identification of square integrable \wifS.
The first one is an extension of Theorem 3 in \cite{OU14a} and the second one is an extension of Theorem 3 in \cite{OU14b} and the second is an extension of the first.

\begin{thm}[\tx{\cite{KH17b}}(2018)] 
Let $(e_n)_{n\in\mab{N}}$ be a CONS of $L^2[0,1]$, where each $e_n(t)$ is bounded.
Then, $a\in\L_1^{1,2}$ is identified from the system $(\sfc{e_n}{a\,\de B})_{n\in\mab{N}}$ of its \sfcss with respect to $(e_n)_{n\in\mab{N}}$.
\label{thm: 12}
\end{thm}

\begin{thm}[\tx{\cite{KH17b}}(2018)] 
Let $(e_n)_{n\in\mab{N}}$ be a CONS of $L^2[0,1]$, where each $e_n(t)$ is bounded and $a\in\L_1^{1,2}$ and $b\in\L_1^{0,2}$.
Define $\de X_t:=a(t)\,\de B_t+b(t)\,dt$.
Then, $a(t)$ and $b(t)$ are identified from the system $(\sfc{e_n}{\de X})_{n\in\mab{N}}$ of \sfcss  of $\de X$ with respect to $(e_n)_{n\in\mab{N}}$.
\label{thm: 21}
\end{thm}


\section{Preliminaries}
\subsection{Notation and terminology}
\label{sub:011}
Throughout this paper, we employ the following setting and terminology:
Let $\s{B}{t}{[0,\infty)}$ be a one-dimensional Brownian motion on a probability space $(\Om,\F,P)$. 
Let $\lambda$ be the Lebesgue measure on $\mab{R}$ and $L$ a constant which satisfies $0<L\le\infty$.
If $L=\infty$ we regard the symbol $[0,L]$ as the infinite interval $[0,\infty)$.
$\L([0,L])$ denotes the \sf of Lebesgue measurable sets on $[0,L]$. 
We say $f:[0,L]\t\Om\to\mab{C}$ is a random function (or measurable stochastic process) on $[0,L]$
if $f$ is $\L([0,L])\ot\F$-measurable.
Note that a random function $f(t)$ is \nC, namely, not necessarily adapted to some filtration for $\s{B}{t}{[0,1]}$, for that matter, not necessarily $\sigma(B_t|t\in [0,\infty))$-measurable for $t\in[0,L]$.
$\ip{f}{g}{}=\ip{f}{g}{L^2([0,L];\mab{C})}$ means the inner product of $f,g\in L^2([0,L];\mab{C})$ defined by $\int_0^L\overline{f}g\,d\lambda$,
where $\overline{f}$ represents the complex conjugate of $f$.
Next, for each $e\in L^2([0,L];\mab{C})$ put $B_t[e]:=\int_0^t e \,dB,\,0\le t\le L$ which is a continuous realization in $t\in[0,L]$.
Also, a real-valued random function whose almost all sample paths are of bounded variation is called a \nc finite variation process or a random function of bounded variation,
and a stochastic process whose almost all paths are (right (resp. left)) continuous is called a (right resp. left) continuous process.
Set $BV[0,L]:=\{v:[0,L]\to\mab{R}\,|\,v \tx{ is of bounded variation }\}$, and
for each $v\in BV[0,L]$ we introduce the following notation:
\begin{enumerate}
\setlength{\itemindent}{-15pt}
\i[$\cdot$] $v_{+}(t)	=\sup\Bigl\{\,\tsum_{j=1}^{n}(v(t_j)-v(t_{j-1}))^{+}\,\Bigl|\,n\in\mab{N}_0, 0=t_0<t_1<\cdots<t_n=t\,\Bigl\}$,\\
\i[$\cdot$] $v_{-}(t)	=\sup\Bigl\{\,\tsum_{j=1}^{n}(v(t_j)-v(t_{j-1}))^{-}\,\Bigl|\,n\in\mab{N}_0, 0=t_0<t_1<\cdots<t_n=t\,\Bigl\}$, 
\i[$\cdot$] \mbox{$v_{\rm tv}(t)	=v_{+}(t)+v_{-}(t)=\tsup\Bigl\{\,\tsum_{j=1}^{n}|v(t_j)-v(t_{j-1})|\,\Bigl|\,n\in\mab{N}_0, 0=t_0<t_1<\cdots<t_n=t\,\Bigl\}$,} 
\end{enumerate}
where $\sup\emptyset:=0,\mab{N}_0=\{0\}\cup\mab{N}$ and $||v||=v_{\rm tv}(L-)$ (the norm of the total variation $\mu_{v_{\rm tv}}$ of the left continuous modification of $v$).
For each measure space $X=(X,\M,\mu)$,
let $L^0(X)=\{\,f:X\to\mab{C}\,|\,f \tx{ is } \M\tx{-measurable and } |f|<\infty \,\,\mu\tx{-{\rm a.e.}}\,\}/\underset{\rm a.e.}{\sim}$,
where $\underset{\rm a.e.}{\sim}$ stands for the equivalence relation on $\{\,f:X\to\mab{C}\,|\,f \tx{ is } \M\tx{-measurable and } |f|<\infty \,\,\mu\tx{-\rm{a.e.}}\,\}$ such that $f\underset{\rm a.e.}{\sim}g\,\,\Leftrightarrow\,\, f(x)=g(x) \,\,\,\mu$-a.a. $x\in X$.
Similarly, we define the space $L^0(X;E)$ of $\M$-measurable functions taking values in a Banach space $E$.
Let $\mab{K}$ be $\mab{R}$ or $\mab{C}$.
In this note, by a CONS we mean an ordered CONS.


\subsection{\og integrals}
\label{sub:00015}

In this subsection,
we give the definitions of two kinds of \og integrals: $\varphi$-integral and {\rm u}-integral.

\begin{dfn}[\og integral] 
Let $f\in L^0(\Om\,;L^2([0,L]\,;\mab{C}))$ and
$(\varphi_m)_{m\in\mathbb{N}}$ be a CONS of $L^2([0,L]\,;\mab{C})$. 
\mpn
\i (\textbf{$\varphi$-integrability})

\noindent We say $f$ is integrable with respect to $\vp$ or $\vp$-integrable if
\begin{equation}
\sum_{m=1}^{\infty}
	\langle \varphi_m,f\rangle_{L^2([0,L])}
	B_L[\varphi_m]
\label{m:2}
\end{equation}
converges in probability.
In that case, \eqref{m:2} is called the \og integral of $f$ with respect to $\vp$ or $\vp$-integral of $f$ and denoted by $\int_{0}^Lf\,d_{\vp}B$. Moreover for each $A\in\L([0,L])$, we say $f$ is $\vp$-integrable on $A$ if $f\ic{A}$ is $\vp$-integrable and then, denote the \og integral by $\int_{A}f\,d_{\vp}B$.

\i (\textbf{{\rm u}-integrability})

\noindent We say $f$ is universally integrable or u-integrable for $L^2([0,L]\,;\mab{K})$ if $f$ is integrable with respect to any CONS $\vp$ in $L^2([0,L]\,;\mab{K})$ and the \og integral $\int_{0}^Lf\,d_{\vp}B$ is independent of the particular choice of $\vp$.
When $f$ is u-integrable for $L^2([0,L]\,;\mab{K})$, we refer to $\int_{0}^Lf\,d_{\vp}B$ as the (universal) \og integral (u-integral) of $f$ for $L^2([0,L]\,;\mab{K})$ and denote it by  $\int_{0}^Lf\,d_{\rm{u}}B$. Moreover for each $A\in\L([0,L])$, we say $f$ is u-integrable on $A$ if $f\ic{A}$ is u-integrable and then, denote the \og integral by $\int_{A}f\,d_{\rm u}B$.
\label{def:0001}
\end{dfn}
\rem\!\! There is a
difference between u-integrability for $L^2([0,1];\mab{R})$ and that for $L^2([0,1];\mab{C})$. For example, by Proposition 6 in \og \cite{Og84a} $\s{B}{t}{[0,1]}$ is not integrable with respect to some ordered CONS of $L^2([0,1];\mab{C})$ composed of the exponential functions $T_n(t)=e^{2\pi \sqrt{-1}nt}, \,n\in\mab{Z}$,
while $\s{B}{t}{[0,1]}$ is u-integrable for $L^2([0,1]\,;\mab{R})$ by Example 2 or Theorem  1 in \og \cite{Og84b}.

\subsection{\wi chaos and \sk integral}
\label{sus:0002}

In this subsection,
we describe the framework of \wi chaos and the definition of \sk integral.
For $n\in\mab{N}$,
we set $L_{\rm sym}^2([0,L]^n)=\{k\in L^2([0,L]^n)\,|\,k \tx{ is symmetric}\,\}$
and $L_{\rm sym}^2([0,L]^0)=L^2([0,L]^0)=\mab{C}$.
For $n\in\mab{N}_0=\{0\}\cup\mab{N}$
we denote by $I_n(k)=\int_{[0,L]^n}k(t_1,...,t_n)\,dB_{t_1}\cdots dB_{t_n}$ the multiple \wii integral of order $n$ of $k\in L^2([0,L]^n)$.
Besides for $i\in\mab{N}_0$, we denote by $L^2_{B}([0,L]^i\t\Om)$ the closed subspace $L^2([0,L]^i\t\Om,\L([0,L])^{\ot i}$
$\ot\F^B,\lambda^{\ot i}\ot P)$ of $L^2([0,L]^i\t\Om)=L^2([0,L]^i\t\Om,\L([0,L])^{\ot i}\ot\F,\lambda^{\ot i}\ot P)$,
where $\F^B$ is the $\sigma$-field generated by $(B_t)_{t\in[0,L]}$.
Then it is easy to see that we can define the multiple \wii integral of order $n$ with $i$-tuple parameters $I_n^{(i)}:L^2([0,L]^{i+n})\to L_B^2([0,L]^i\times\Omega)$ such that $I_n^{(i)}\Bigl(\,\textstyle\sum\limits_{j=1}^{r}e_j\otimes k_j\,\Bigl)=\textstyle\sum\limits_{j=1}^{r}e_j\otimes I_n(k_j) \text{ for any } r\in\mathbb{N},
e_j\in L^2([0,L]^i)\text{ and }k_j\in L^2([0,L]^n)$. Denote $I_n^{(1)}$ by $\widehat{I}_n$, in particular.
Then, the \wii Theorem implies that
$L_B^2([0,L]^i\times\Omega)$ is decomposed into the orthogonal direct sum as follows:
\begin{equation*}
L_B^2([0,L]^i\times\Omega)=\displaystyle\bigoplus_{n=0}^{\infty}I^{(i)}_n(L^2_{\cdot;\rm sym}([0,L]^{i+n})),
\end{equation*}
where $L^2_{\cdot;\rm sym}([0,L]^{i+n})=\{k\in L^2([0,L]^{i+n})\,\,|\,\,k(t_1,...,t_i;\cdot)\in L^2_{\rm sym}([0,L]^n) \,\,\,  \text{for a.a. } (t_1,...,t_i)$
$\in [0,L]^{i}\}\,\tx{ and }L^2_{\cdot;\rm sym}([0,L]^{0+0})=\mathbb{C}$.
Thus, any \wi functional $f\in L_B^2([0,L]^i\times\Omega)$ is uniquely expressed in $L_B^2([0,L]^i\times\Omega)$ as
\begin{equation}
f(t)
=\sum_{n=0}^{\infty}I_n^{(i)}(k_n^f(t;\cdot)),
\label{a:0011}
\end{equation}
where $k_n^f(\cdot;\cdot)\in L^2_{\cdot;\rm sym}([0,L]^{i+n})$ is called the kernel of order $n$ for $f$, and $k_0^f(t) = E(f(t))$.

\noindent In addition, \eqref{a:0011} is  rephrased as
$
f(t)=\tsum_{n=0}^{\infty}I_n(k_n^f(t;\cdot)) \,\, \text{in}\,\, L^2(\Omega) \,\,\text{ for almost all } t \in [0,L]^i
$
and the following commutativity between the inner product on $L^2([0,L]^i)$ and the summation of the expansion holds:
$$
\langle f,h\rangle_{L^2([0,L]^i)}=\displaystyle\sum_{n=0}^{\infty}\langle f_n,h\rangle_{L^2([0,L]^i)} \quad \text{in}\,\, L^2(\Omega)\quad\text{ for all }h\in L^2([0,L]^i).
$$

\begin{dfn}[\sosS]
For each $i\in\mab{N}_0$ and $r\in[0,\infty)$, we define the \sos $\L_i^{r,2}$ by
\begin{equation*}
\L_i^{r,2}
=\Bigl\{\,
	f\in L_B^2([0,L]^i\t\Om)	\,\,\Bigl|\,\,
	|f|_{i,r,2}^2:
	=\sum_{n=0}^{\infty}
		(n+1)^rn!\,|k_n^f(\cdot;\cdot)|_{L^2[0,L]^{i+n}}^2
	<\infty
  \,\Bigl\}.
\end{equation*}
For each $r\in(-\infty,0)$, we define the norm on $L_B^2([0,L]^i\t\Om)$ by
\begin{equation*}
	|f|_{i,r,2}^2:
	=\sum_{n=0}^{\infty}
		(n+1)^rn!\,|k_n^f(\cdot;\cdot)|_{L^2[0,L]^{i+n}}^2
\end{equation*}
and define the \sos $\L_i^{r,2}$ as the completion of $L_B^2([0,L]^i\t\Om)$ with respect to the norm $|\cdot|_{i,r,2}$.
\label{def:00111}
\end{dfn}
Then, the following hold:
\begin{enumerate}
\i[$\cdot$] $\L_i^{r,2}$ becomes a Hilbert space with respect to the norm $|\cdot|_{i,r,2}$.
\i[$\cdot$] $\L_i^{0,2}=L_B^2([0,L]^i\t\Om)$.
\i[$\cdot$] \mbox{$\L_i^{r,2}$ monotonically decreases with respect to $r:r_1<r_2 \text{ implies } \L_i^{r_1,2} \supset\L_i^{r_2,2}$.}
\end{enumerate}

\begin{dfn}[$H$-derivative]\bt
Let $i\in\mathbb{N}_0,\,r\in\mab{R}$ and $f\in\L_i^{r,2}$ which
is represented as follows:
$
f(s)=\tsum_{n=0}^{\infty}I_{n}\bigl(\,k_n^f(s;\cdot)\,\bigl) \,\,\tx{ in }\L_i^{r,2}.
$
Then,
\begin{equation}
\displaystyle\sum_{n=1}^{\infty}nI_{n-1}\bigl(\,k_n^f(s;\cdot,t)\,\bigl)
\label{s:53}
\end{equation}
converges in $\mathcal{L}_{i+1}^{r-1,2}$. We denote \eqref{s:53} by $D_tf(s)$.
\label{z:90}
\end{dfn}

\begin{dfn}[\sk integral] 
Let $i\in\mathbb{N},\, r\in\mab{R}$ and $f\in\L_i^{r,2}$ which
is represented as follows: 
$
f(t_1,...,t_i)
=$
$
\tsum_{n=0}^{\infty}
I_n\bigl(\,
	k_n^f(t_1,...,t_i\,;\, \dot{s}_1,...,\dot{s}_n)
\,\bigl)\,\,\tx{ in }\L_i^{r,2}.
$
Then, for $j\in\{1,\ldots,i\}$
\begin{equation}
\sum_{n=0}^{\infty}
I_{n+1}\bigl(\,
	k_n^f(
		t_1,...,\overset{j\text{-th}}{\overset{\vee}{\dot{s}_{n+1}}},...,t_i
		\,;\, 
		\dot{s}_1,...,\dot{s}_n
	)
\,\bigl)
\label{s:48}
\end{equation}
converges in $\mathcal{L}_{i-1}^{r-1,2}$. We call \eqref{s:48} the \sk integral of $f(t_1,...,t_i)$ with respect to $t_j$,
and denote it by $\int_{0}^Lf(t_1,...,t_i)\,\de B_{t_j}$.
\label{dfn:92}
\end{dfn}

\begin{dfn}[Hilbert-Schmidt integral operators] 
Given $K\in L^2([0,L]^2)$ and $r\in[0,\infty)$,
we define the operators $T_K$ on $\d\sqcu_{n=0}^{\infty}L^2([0,L]^{1+n})$ and $\widehat{T}_K$ on
$\L_1^{r,2}$ by
$$
T_Kk_n(t;\cdot)=(T_Kk_n(\cdot;\cdot))(t)=\int_{0}^LK(t,s)k_n(s;\cdot)\,ds,
$$
$$
\widehat{T}_Kf(t)=(\widehat{T}_Kf)(t)=\int_{0}^LK(t,s)f(s)\,ds.
$$
\end{dfn}
\rem  For notational simplicity,
we use the same symbol $T_K$ to represent these two operators since one can distinguish between them in the context.


\section{Representations of \og integrals}
\label{ssec:3.1001}

In this section, we give sufficient conditions for random functions to be \og integrable and explicit representations of the integrals, which are necessary for proving the main theorems in Sections 6 and 7.

\subsection{Integration by parts}

To begin with, the next theorem follows the It\^o-\ni theorem (Theorem 4.1 in \cite{IN}).
One can prove it by the same argument as in the proof of Theorem 5.1 in \cite{IN}.

\begin{thm}[It\^o-\nI] 	
Let $(\vp_m)_{m\in\mab{N}}$ be a CONS of $L^2([0,L]\,;\mab{R})$. 
For $e\in L^2([0,L]\,;\mab{C})$
and $t\in[0,L]$,
the following holds:

\begin{align*}
\sum_{m=1}^{\infty}B_L[\vp_m]\ip{e\mathsf{1}_{[0,s]}}{\vp_m}{}=B_s[e] \q\text{uniformly in $s\in[0,t]$ almost surely.}
\end{align*}
\label{thm: 3.101}
\end{thm}

\begin{prf}
From linearities we can assume $K(u,s)=e(u)\ic{u\le s}$ is real-valued. Set $(X_m)_{m\in\mathbb{N}}$ as a sequence of $C([0,t]\,;\mab{R})$-valued random variables
such that\\
$X_m(s)=\ip{K(\cdot,s)}{\vp_m}{}B_L[\vp_m]$,
and set $X$ as a $C([0,t]\,;\mab{R})$-valued random variable
such that
$X(s):=\int_0^LK(u,s)\,dB_u$.
Let $S_M=\tsum_{m=1}^MX_m$ for each $M\in\mab{N}$.
We are to show $\tlim_{M\to\infty} S_M=X \,\, \text{almost surely.}$
Obviously, $(X_m)_{m\in\mathbb{N}}$ satisfies
\begin{enumerate}
\i[(a)] $(X_m)_{m\in\mathbb{N}}$ is independent,
\i[(b)] each $X_m$ is symmetrically distributed,
\end{enumerate}
so that $\tlim_{M\to\infty} S_M=X \,\, \text{almost surely}$
if and only if
$\,\a z\in C([0,t])^*\,\tlim_{M\to\infty}{}_{C([0,t])^*}{\ip{z}{S_M}{C([0,t])}}={}_{C([0,t])^*}{\ip{z}{X}{C([0,t])}}$ in probability by Theorem 4.1 in \cite{IN}.
Fix any $z\in C([0,t])^*$, by the Riesz-Markov theorem  there exists a regular \bo finite signed measure $\mu_{z}$ such that $\a f\in C([0,t])\,\,\blf{C([0,t])^*}{z}{f}{C([0,t])}=\int_{[0,t]}f\,d\mu_{z}$. Then we have 
\begin{align*}
\aex{\blf{C([0,t])^*}{z}{X-S_M}{C([0,t])}}
=&\AEx{\int_{[0,t]}(X-S_M)(s)\,d\mu_{z}(s)}\\
\le&\EX{\int_{[0,t]}|X-S_M|(s)\,d|\mu_{z}|(s)}\\
=&\int_{[0,t]}E(\,|X-S_M|(s)\,)\,d|\mu_{z}|(s),
\end{align*}
and
\begin{align*}
E(\,|X-S_M|(s)\,)
\le&E(\,|X-S_M|^2(s)\,)^{\frac{1}{2}}\\
=&\Bigl(\,\lim_{L\to\infty}E\Bigl|\!\!\!\sum_{\,\,\,m=M+1}^{L}\ip{K(\cdot,s)}{\vp_m}{}B_L[\vp_m]\,\Bigl|^2\,\Bigl)^{\frac{1}{2}}\\
=&\Bigl(\!\!\!\sum_{\,\,\,m=M+1}^{\infty}
				\ip{K(\cdot,s)}{\vp_m}{}^2\,\Bigl)^{\frac{1}{2}}\\
&\!\!\!\!\!\begin{cases}
\underset{M\to\infty}{\longrightarrow} 0 & ,\tx{a.a. }\, s\in [0,t]\\
\,\,\le\q |K(\cdot,s)|_{L^2([0,L])}\le(\int_0^te(u)^2\,du)^{\frac{1}{2}}<\infty& ,\tx{a.a. }\,  s\in [0,t].\\
\end{cases}
\end{align*}
Because $|\mu_z|$ is finite, we can apply the bounded convergence theorem and we gain
$$
\lim_{M\to\infty}\aex{\blf{C([0,t])^*}{z}{X-S_M}{C([0,t])}}=0.
$$
Therefore, we conclude 
$$
\lim_{M\to\infty}\blf{C([0,t])^*}{z}{S_M}{C([0,t])}=\blf{C([0,t])^*}{z}{X}{C([0,t])} \q \text{in probability.}
$$
\end{prf}

Now, in light of \tr{thm: 3.101} we have the next integration by parts lemma.
\begin{lem} 
Assume $e\in L^2([0,L]\,;\mab{C})$ has compact support. 
Let $a(t)$ be a \nc finite variation process on $[0,L]$ and $0\le s<t\le L$.
Then, $ae$ is u-integrable on $[s,t]$ for $L^2([0,L]\,;\mab{R})$ and the \og integral converges almost surely and it is given by 
\begin{align}
\int_s^{t} ae\,d_{\rm u}B=a(t-)B_t[e]-a(s+)B_s[e]-\int_{(s,t)} B_u[e]\,da(u).
\label{al:0.77}
\end{align}
\label{lem:0.7}
\end{lem}
\begin{prf}\bt
Fix any CONS $(\vp_m)_{m\in\mab{N}}$ of $L^2([0,L]\,;\mab{R})$.
For now, the symbol of equality '=' means both sides of the expression are equal almost surely.

Case 1, the case of $\,t<\infty$ :
By \leB's integration by parts formula,  for $M\in\mab{N}$ we have
\begin{align}
&\sum_{m=1}^{M} B_L[\vp_m] \int_s^{t} ae\vp_m \,d\lambda\nonumber\\
=&\sum_{m=1}^{M}B_L[\vp_m] \Bigl(\,a(t-) \ip{e\ic{[0,t]}}{\vp_m}{}-a(s+) \ip{e\ic{[0,s]}}{\vp_m}{}-\int_{(s,t)} \ip{e\ic{[0,u]}}{\vp_m}{}
\,da(u)\,\Bigl)\nonumber\\
=&a(t-)\sum_{m=1}^{M}B_L[\vp_m]\ip{e\ic{[0,t]}}{\vp_m}{}
-a(s+)\sum_{m=1}^{M}B_L[\vp_m]\ip{e\ic{[0,s]}}{\vp_m}{}\nonumber\\
&-\sum_{m=1}^{M}B_L[\vp_m]\int_{(s,t)}\ip{e\ic{[0,u]}}{\vp_m}{}\,da(u),
\label{al:0.888}
\end{align}
the series of the first and second terms of \eqref{al:0.888} are \wi expansions, so that they converge almost surely, and the third term of \eqref{al:0.888} converges by \tr{thm: 3.101} to $-\int_{(s,t)} B_u[e]\,da(u)$ almost surely. Therefore we get
$$
\sum_{m=1}^{\infty} B_L[\vp_m] \int_s^{t} ae\vp_m \,d\lambda=a(t-)B_t(e)-a(s+)B_s(e)-\int_{(s,t)} B_u[e]\,da(u) \q\tx{almost surely.}
$$
\\
Case 2, the case of $\,t=\infty$ : 
Because $e$ has compact support, there exists $L_0\in\mab{N}$ such that $\a u\ge L_0 \,\, e(u)=0$.
By the definition of \og integral and Case 1 we have
\begin{align*}
\int_s^{\infty} ae\,d_{\rm u}B
=&\int_s^{L_0} ae\,d_{\rm u}B
=a(L_0-)B_{L_0}[e]-a(s+)B_s[e]-\int_{(s,L_0)} B_u[e]\,da(u).
\end{align*}
Then we have \eqref{al:0.77} since
\begin{align*}
\int_{[L_0,\infty)} B_u[e]\,da(u)
=\int_{[L_0,\infty)} B_{L_0}[e]\,da(u)
=&B_{L_0}[e]\,(a(\infty)-a(L_0-)).
\end{align*}
\end{prf}

\subsection{Integrals of \wifs}
Next, we give some propositions in the case that the integrand is \wi functional.
The next theorem is a known result, Proposition 6 in \og \cite{Og84a}, which gives a necessary and sufficient condition for  $f\in\L_1^{1,2}$ to be $\vp$-integrable and a representation of the integral.
\begin{thm}[Ogawa]\bt
Let $f\in\L_1^{1,2}$ and $\vp$ be a CONS of $L^2([0,1]\,;\mab{C})$.
The following are equivalent:
\begin{enumerate}	
\i[$(i)$] $f$ is $\vp$-integrable.
\i[$(ii)$] $\displaystyle\sum_{m=1}^{\infty}\int_0^1\int_0^1D_tf(s)\vp_m(t)\overline{\vp_m(s)}\,dt\,ds\,$ converges in probability. 
\end{enumerate}
Moreover, in this case the \og integral is given by
$$
\int_0^1f(t)\,d_{\vp}B_t=\int_0^1f(t)\,\de B_t+\displaystyle\sum_{m=1}^{\infty}\int_0^1\int_0^1D_tf(s)\vp_m(t)\overline{\vp_m(s)}\,dt\,ds.
$$
\label{thm:3.10}
\end{thm}

From Theorem \ref{thm:3.10}, the next assertion holds.
\begin{cor}\bt
Let $f\in\L_1^{1,2}$ and suppose $Df\in L^2([0,L]\,;\C)$ is of trace class almost surely, then
$f$ is u-integrable for $L^2([0,L];\mab{C})$ and the \og integral is given by
\begin{align}
\int_0^Lf(t)\,d_{\rm u}B_t=\int_0^Lf(t)\,\de B_t+\tx{tr}(Df).
\label{al:3.4}
\end{align}
\label{pr:3.51}
\end{cor}
\rem If $K\in L^2([0,L]^2) \tx{ and }g\in\L_1^{1,2}$, then $f:=T_Kg\in\L_1^{1,2}$,
$D_tf(s)=\ip{\,\overline{D_tg(\cdot)}}{K(s,\cdot)\,}{}$, and $T_{Df}=T_{Dg}T_{K^*}$, which is of trace class.
Therefore the universal \og integral of $f$ is given by \eqref{al:3.4}. But we can't see directly from this fact the explicit expression of the trace of $Df$. We have the following in this context.

The next theorem says an $\L_1^{1,2}$ class \wi functional is u-integrable for $L^2([0,L];\mab{C})$ when its kernels in \wii decomposition have integral representations with Hilbert-Schmidt kernels.

\begin{prp}
Let $(k_n)_{n\in\mab{N}_0}\in\tprod_{n=0}^{\infty}L^2([0,L]^{1+n})$ and $(K_n)_{n\in\mab{N}_0}$ be a sequence of $L^2([0,L]^2)$.
Suppose the following holds:
\begin{align}
\sum_{n=0}^{\infty}(n+1)!|\hat{k}_n|_{L^2([0,L]^{1+n})}^2|K_n|_{L^2([0,L]^2)}^2<\infty.
\label{al:10.9}
\end{align}
Then $F:=\tsum_{n=0}^{\infty}\widehat{I}_nT_{K_n}k_n=\tsum_{n=0}^{\infty}T_{K_n}\widehat{I}_nk_n\in\L_1^{1,2}$ is u-integrable for $L^2([0,L]\,;\mab{C})$ and the \og integral converges in $L^2(\Om)$ and it is given  by 
\begin{align*}
\int_{0}^LF(t)\,d_{\rm u}B_t
&=\int_{0}^LF(t)\,\de B_t+\sum_{n=1}^{\infty}nI_{n-1}(\hat{k_n}^*\otimes_2K_n)\\
&=\int_{0}^LF(t)\,\de B_t+\sum_{n=1}^{\infty}\int_{[0,L]^2}K_n(t,s)D_tI_n((\hat{k}_n(s;\cdot))\,ds\,dt,
\end{align*}
where $\hat{k}_n(t;\cdot)=\tx{Sym}_n(k_n(t;\cdot))$ and $\hat{k_n}^*\in L^2([0,L]^{n+1})$ defined by $\hat{k_n}^*(\cdot;t)=\hat{k_n}(t;\cdot)$.
\label{thm: 11}
\end{prp}
\rem
If for some $r\in\mab{R}$, 
$
a=\tsum_{n=0}^{\infty}I_n(k_n)\in\L_1^{r,2}\,\tx{ and }\,\tsup_{n\in\mab{N}_0}n^{\frac{1-r}{2}}|K_n|_{L^2([0,L]^2)}<\infty,
$
then (\ref{al:10.9}) holds.

\begin{prf} 
Note that $\widehat{I}_nT_{K_n}\hat{k}_n=\widehat{I}_nT_{K_n}k_n$ because $\widehat{I}_nT_{K_n}=T_{K_n}\widehat{I}_n$.
We have
\begin{align*}
|F|_{1,1,2}^2
&=\sum_{n=0}^{\infty}(n+1)!\int_{[0,L]^{n+1}}|\ip{\overline{K_n(t,\cdot)}}{\hat{k}_n(\cdot\,;\,u_1,...,u_n)}{}|^2\,dtdu_1...du_n\\
&\le\sum_{n=0}^{\infty}(n+1)!|\hat{k}_n|_{L^2([0,L]^{1+n})}^2|K_n|_{L^2([0,L]^2)}^2<\infty
\end{align*}
by the \scI, thus $F\in\L_1^{1,2}$.
Fix any CONS $(\vp_m)_{m\in\mab{N}}$ of $L^2([0,L]\,;\mab{C})$.
By Theorem \ref{thm:3.10}, we should investigate the term $\tsum_{n=1}^{\infty}nI_{n-1}(C_n^M)$, where $M\in\mab{N}$ and
\begin{align*}
C_n^M(u_1,...,u_{n-1})
&=\sum_{m=1}^{M}\int_{0}^L\int_{0}^Lk_{n}^F(t\,;\,u_1,...,u_n)\overline{\vp_m(t)}\vp_m(u_{n})\,dt\,du_n\\
&=\sum_{m=1}^{M}
	\int_{0}^L\int_{0}^L\int_{0}^LK_n(t,s)\hat{k}_{n}(s\,;\,u_1,...,u_n)\overline{\vp_m(t)}\vp_m(u_{n})\,ds\,dt\,du_n.
\end{align*}
Because 
\begin{align*}
&\int_{[0,L]^3}|K_n(t,s)\hat{k}_n(s\,;\,u_1,...,u_n)\overline{\vp_m(t)}\vp_m(u_n)|\,dsdtdu_n\\
\le&\Bigl(\,\int_{[0,L]^3}|K_n(t,s)\vp_m(u_n)|^2\,dsdtdu_n\int_{[0,L]^3}|\hat{k}_n(s\,;\,u_1,...,u_n)\overline{\vp_m(t)}|^2\,dtdsdu_n\,\Bigl)^{\frac{1}{2}}\\
=&\Bigl(\,|K_n|_{L^2([0,L]^2)}^2\int_{[0,L]^2}|\hat{k}_n(s\,;\,u_1,...,u_n)|^2\,dsdu_n\,\Bigl)^{\frac{1}{2}}<\infty
\end{align*}
for almost all $\,(u_1,...,u_{n-1})\in [0,L]^{n-1}$, by Fubini's theorem we have
\begin{align*}
C_n^M(u_1,...,u_{n-1})&=\sum_{m=1}^{M}\int_{0}^L\int_{0}^L\ip{K_n(\cdot,s)}{\vp_m}{}\hat{k}_{n}(s\,;\,u_1,...,u_n)\vp_m(u_{n})\,du_nds.
\end{align*}
Now, we are to show $\tsum_{n=1}^{\infty}nI_{n-1}(C_n^M)$ converges to
$
\tsum_{n=1}^{\infty}nI_{n-1}(\hat{k_n}^*\otimes_2K_n)
$
in $L^2(\Om)$ as $M\to\infty$.
First, we have $\tsum_{n=1}^{\infty}nI_{n-1}\Bigl(\,\hat{k_n}^*\otimes_2K_n\,\Bigl)\in L^2(\Om)$
since
\begin{align*}
&\Bigl|\,\sum_{n=1}^{\infty}nI_{n-1}(\hat{k_n}^*\otimes_2K_n)\,\Bigl|_{0,0,2}^2\\
=&\sum_{n=1}^{\infty}n^2(n-1)!\int_{[0,L]^{n-1}}\Bigl|\,\int_{[0,L]^2}\hat{k}_n(s;u_1,...,u_n)K_n(u_n,s)\,du_nds\,\Bigl|^2du_1...u_{n-1}\\
\le&\sum_{n=1}^{\infty}nn!|\hat{k}_n|_{L^2([0,L]^{1+n})}^2|K_n|_{L^2([0,L]^2)}^2<\infty
\end{align*}
by the assumption (\ref{al:10.9}) and the \scI.
Next, we have
\begin{align*}
&\Bigl|\,\sum_{n=1}^{\infty}nI_{n-1}(\hat{k_n}^*\otimes_2K_n-C_n^M)\,\Bigl|_{0,0,2}^2\\
=&\sum_{n=1}^{\infty}nn!\!\int_{[0,L]^{n-1}}\Bigl|\int_{[0,L]^2}\!\!\hat{k}_n(s;u_1,...,u_n)\Bigl(K_n(u_n,s)\\
&\q-\sum_{m=1}^{M}\ip{K_n(\cdot,s)}{\vp_m}{}\vp_m(u_n)\Bigl)\,du_nds\Bigl|^2\!du_1...u_{n-1}\\
\le&\sum_{n=1}^{\infty}nn!|\hat{k}_n|_{L^2([0,L]^{1+n})}^2\int_{0}^L\,\Bigl|\,K_n(\dot{u}_n,s)-\sum_{m=1}^{M}\ip{K_n(\cdot,s)}{\vp_m}{}\vp_m(\dot{u}_n)\,\Bigl|_{L^2([0,L],\,du_n)}^2\,ds\\
=&\sum_{n=1}^{\infty}nn!|\hat{k}_n|_{L^2([0,L]^{1+n})}^2\int_{0}^L\Bigl(\,|K_n(\cdot,s)|_{L^2([0,L])}^2-\sum_{m=1}^{M}|\ip{K_n(\cdot,s)}{\vp_m}{}|^2\,\Bigl)\,ds
\end{align*}
and by the monotone convergence theorem, we have 
$$
\int_{0}^L\Bigl(\,|K_n(\cdot,s)|_{L^2([0,L])}^2-\sum_{m=1}^{M}|\ip{K_n(\cdot,s)}{\vp_m}{}|^2\,\Bigl)\,ds\underset{M\to\infty}{\searrow}0
$$
for each $n\in\mab{N}$. By (\ref{al:10.9}) and applying the monotone convergence theorem again, we get

\noindent
$
\tlim_{M\to\infty}\Bigl|\,\tsum_{n=1}^{\infty}nI_{n-1}(\hat{k_n}^*\otimes_2K_n-C_n^M)\,\Bigl|_{0,0,2}^2=0.
$
Therefore, $\tlim_{M\to\infty}\tsum_{n=1}^{\infty}nI_{n-1}(C_n^M)=$

\noindent
$\tsum_{n=1}^{\infty}nI_{n-1}(\hat{k_n}^*\otimes_2K_n)$
in $L^2(\Om)$. Finally we have
\begin{align*}
\sum_{n=1}^{\infty}nI_{n-1}(\hat{k_n}^*\otimes_2K_n)
&=\sum_{n=1}^{\infty}\int_{[0,L]^2}K_n(t,s)nI_{n-1}(\hat{k}_n(s;\cdot,t))\,dtds\\
&=\sum_{n=1}^{\infty}\int_{[0,L]^2}K_n(t,s)D_tI_n((\hat{k}_n(s;\cdot))\,dtds,
\end{align*}
which completes the proof of Proposition \ref{thm: 11}.
\end{prf}

\begin{cor}
Suppose $F$ is represented by 
\begin{align*}
F(t)
=T_Kf(t)
=\int_{0}^LK(t,s)f(s)\,ds,
\end{align*}
where $K\in L^2([0,L]^2)$ and $f\in\L_1^{1,2}$. Then $F$ is u-integrable for $L^2([0,L]\,;\mab{C})$ and the \og integral converges in $L^2(\Om)$ and it is given by
\begin{align*}
\int_{0}^LF(t)\,d_{\rm u}B_t
&=\int_{0}^LF(t)\,\de B_t+\int_{0}^L\int_{0}^LK(t,s)D_tf(s)\,ds\,dt.
\end{align*}
\label{cor:13.5}
\end{cor}
\rem We can prove Corollary \r{cor:13.5} directly by the same argument as
that of the proof of Proposition \r{thm: 11}.

From Corollary \r{cor:13.5}, we have the next statement, which says the \og integral of product of a locally absolutely continuous \wi functional and a deterministic square integrable function with compact support.

\begin{cor} 
Let $a(t)$ be a real valued random function on $[0,L]$ and $e\in L^2([0,L])$ has compact support.
Suppose $a(t)$ satisfies the following:
\begin{enumerate}
\i $a(t)$ is locally absolutely continuous in $t\in[0,L]$ almost surely.
\i $a'(t):=\frac{d}{dt}a(t)\in\mathcal{L}_1^{1,2}$ and $a(0)\in\mathcal{L}_0^{1,2}$.
\end{enumerate}
Then, $ae$ is u-integrable for $L^2([0,L]\,;\mab{C})$ and the \og integral converges in $L^2(\Om)$ and it is given by
\begin{align}
\int_0^Lae(t)\,d_{\rm u}B_t=\int_0^Lae(t)\,\de B_t+\int_0^{L}\Bigl(\,\int_0^{t}D_ta'(s)\,ds+D_ta(0)\,\Bigl)e(t)\,dt.
\label{al:14.1}
\end{align}
\label{cor:13.7}
\end{cor}
\begin{prf} 
Because $a(t)$ is represented by
$
a(t)=\int_0^ta'(s)\,ds+a(0),
$
set
\[
K_1(t,s)=e(t)\ic{s\le t}\,\,\tx{and }\, K_2(t,s)=
\begin{cases}
\frac{1}{L}e(t)&,\tx{if } L<\infty\\
e(t)\ic{s\le 1}&,\tx{if } L=\infty
\end{cases}
\]
and then
$K_{2}\in L^2([0,L]^2)$ and $K_{1}\in L^2([0,L]^2)$, for $e$ has compact support.
Then, by Corollary \r{cor:13.5}, 
$ae=T_{K_{1}}a'+T_{K_{2}}a(0)$ is u-integrable for $L^2([0,L]\,;\mab{C})$ and \eqref{al:14.1} holds.
\end{prf}
\rem If $L<\infty$, then $a(t)$ is a \nc finite variation process and it has already proved that $ae$ is u-integrable for $L^2([0,L]\,;\mab{R})$ by Lemma \r{lem:0.7}. This corollary gives another representation as \eqref{al:14.1} 
and shows the left hand side of \eqref{al:14.1} converges in $L^2(\Om)$ as well as almost surely to the right hand side.

On the other hand, from Corollary \r{cor:13.5} and Theorems 5.1,5.2 and 5.3 in \cite{KH17a}, the following statement about the \og integrability holds.

\begin{thm}\bt
Let $e:[0,1]\to\mab{R}$ and a CONS $(\vp_m)_{\in\mab{N}}$ of $L^2([0,1];\mab{R})$ and a process
\begin{align}
a(t):=\int_0^tf(s)\,\de B_s+T_Kh(t)+a(0)\,\,,t\in[0,1],
\label{eq:17}
\end{align}
where $f\in\L_1^{2,2}, h\in\L_1^{1,2}, a(0)\in\L_0^{1,2}$ and $K\in L^2[0,1]^2$.
Let  (C.1),(C.2) and (C.3) be conditions on $e(t)$, $(\vp_m)_{\in\mab{N}}$ and $f(t)$ stated in Subsection \r{sub: 1.1}.
Then, we have the following assertions (A),(B) and (C):
\begin{enumerate}
\i[(A)] \, If (C.1) or (C.2) hold, then $ea$ is $\varphi$-integrable.
\i[(B)] \, If (C.3) holds, then $ea$ is u-integrable for $L^2([0,1];\mab{R})$.
\i[(C)] \, If the assumption of (A) or (B) holds, the integral of $ea$ converges in $L^2(\Om)$ and it is given by 
$$
\int_0^1e(t)a(t)\,\de B_t+\frac{1}{2}\int_0^1e(t)f(t)\,dt
$$
\begin{align}
+\int_0^1e(t)\Bigl(\,\int_0^tD_tf(s)\,\de B_s+\int_0^1K(t,s)D_th(s)\,ds+D_ta(0) \,\Bigl)\,dt.
\label{al:17:01}
\end{align}
In particular, when $K(t,s)=\ic{s\le t}$, $a(t)$ defines an S-type It\^o process
\begin{align}
a(t)=\int_0^tf(s)\,\de B_s+\int_0^th(s)\,ds+a(0)\,,t\in[0,1]
\label{al:15.7}
\end{align}
and \eqref{al:17:01} is rephrased as 
$$
\int_0^1e(t)a(t)\,\de B_t+\frac{1}{2}\int_0^1e(t)f(t)\,dt
+\int_0^1e(t)\Bigl(\,\int_0^tD_t\,\de a(s)+D_ta(0) \Bigl)\,dt,
$$ 
where $\int_0^tD_t\,\de a(s)$ denotes $\int_0^tD_t\,f(s)\,\de B_s+\int_0^tD_t\,h(s)\,ds$.
\end{enumerate}
\label{thm: 17.02}
\end{thm}


\section{\sfcis and identification problems}
\label{sec:4}

In this section, 
we state the definitions of \sfcis introduced by \og and introduce the definitions of identification of random functions.
From now to the end of this paper (except for Appendixes), let $(e_n)_{n\in\mab{N}}$ be a fixed CONS of $L^2([0,L]\,;\mab{C})$,
$a\in L^0([0,L]\t\Om\,;\C)$ and $b\in L^0(\Omega\,;L^2([0,L]\,;\mab{C}))$.


\subsection{Definitions of \sfciS}

We state the definitions of \sfcis of a stochastic differential by $a(t)$ and $b(t)$ with respect to $(e_n)_{n\in\mab{N}}$.

\begin{dfn}[SFC-$\rm{\textsc{O}}_{\bf\vp}$  of stochastic differential] 
Let $(\vp_m)_{m\in\mab{N}}$ be a CONS of $L^2([0,L]\,;\,\mab{C})$.
Suppose $\overline{e_n}a$ is $\vp$-integrable for every $n\in\mathbb{N}$. 
We define the $n$-th \sfci of $\vp$-\og type (\sfcoP) $\sfc{e_n}{\dpy}$ of the stochastic differential 
$
\dpy_t=a(t)\,d_{\vp}B_t+b(t)\,dt\,\,,t\in [0,L]
$
with respect to $(e_n)_{n\in\mathbb{N}}$ by
$$
\sfc{e_n}{\dpy}:=\int_{0}^L\overline{e_n(t)}\,\dpy_t=\int_{0}^L\overline{e_n(t)}a(t)\,d_{\vp}B_t+\int_{0}^L\overline{e_n(t)}b(t)\,dt.
$$
In Particular, in the case of $b=0$, $\sfc{e_n}{\dpy}=\sfc{e_n}{a\,d_{\vp}B}$ is also called the \sfcop of $a(t)$.
\end{dfn}

\begin{dfn}[SFC-$\rm{\textsc{O}}_{\bf u}$ of stochastic differential] 
Suppose $\overline{e_n}a$ is u-integrable for every $n\in\mathbb{N}$. 
We define the $n$-th \sfci of (universal) \og type (\sfcoU) $\sfc{e_n}{\du Y}$ of the stochastic differential 
$
d_{\rm u}Y_t=a(t)\,d_{\rm u}B_t+b(t)\,dt\,\,,t\in [0,L]
$
with respect to $(e_n)_{n\in\mathbb{N}}$ by
$$
\sfc{e_n}{\du Y}:=\int_{0}^L\overline{e_n(t)}\,d_{\rm u}Y_t=\int_{0}^L\overline{e_n(t)}a(t)\,d_{\rm u}B_t+\int_{0}^L\overline{e_n(t)}b(t)\,dt.
$$
In particular, in the case of $b=0$, $\sfc{e_n}{\du Y}=\sfc{e_n}{a\,\du B}$ is also called the \sfcou of $a(t)$.
\end{dfn}

\begin{dfn}[\sfcs of a stochastic differential]\bt
Suppose any of the following two conditions:
\begin{enumerate}
\i $e_na\in\mathcal{L}_1^{1,2}$ for every $n\in\mathbb{N}$.
\i $e_na\in\mathcal{L}_1^{0,2}$ for every $n\in\mathbb{N}$ and $b\in\mathcal{L}_1^{0,2}$.
\end{enumerate}
We define the $n$-th \sfci of \sk type (\sfcS) $\sfc{e_n}{\de X}$ of the stochastic differential 
$
\de X_t=a(t)\,\de B_t+b(t)\,dt\,\,,t\in [0,L]
$
with respect to $(e_n)_{n\in\mathbb{N}}$ by
$$
\sfc{e_n}{\de X}:=\int_{0}^L\overline{e_n(t)}\,\de X_t=\int_{0}^L\overline{e_n(t)}a(t)\,\de B_t+\int_{0}^L\overline{e_n(t)}b(t)\,dt,
$$
which is an element in $\L_0^{-1,2}$ in the case of 2.
In particular, in the case of $b=0$, $\sfc{e_n}{\de X}=\sfc{e_n}{a\,\de B}$ is also called the \sfcs of $a(t)$.
\end{dfn}

\subsection{B-dependency and Constructiveness on identification}
\label{sub:3.2}

In this subsection, the following notion as the foundation is introduced.
We define some types of identification of a random function.
In this subsection except Example \ref{ex: 4.5}, we use notions introduced in Appendix A.
Hereafter, we forget all the setting given in Subsection 2.1 and we start discussion without any preliminaries: 
Let $(\Om,\F,P)$ be a probability space and $I$ a real interval and set $\pmb{B}_I(\Om,\F,P):=\{\,\s{B}{t}{I}\,|\,\s{B}{t}{I}$ is a Brownian motion on
$(\Om,\F,P)\,\}$.
Let $\L$ be an assigned language of first-order logic defined in Appendix A.

\begin{dfn}
Let $\Lambda$ be an arbitrary set and 
put $\mathscr{L}(I\times\Om)=\{f:I\times\Omega\to\C\,|\,f$
is essentially measurable on
$I\times\Om$
and $f$ is finite $\mu\ot P\tx{-{\rm a.e.}}\,\}$,
which may be divided by some equivalent relation.
We denote by $\mathcal{X}$ the totality of maps from a subset of $\mathscr{L}(I\times\Om)$ to $L^0(\Om)^{\Lambda}$, where $L^0(\Omega)^{\Lambda}$ denotes the family of sequences indexed by $\Lambda$ of random variables on $(\Om,\F,P)$. 
Let $X:\pmb{B}_I(\Om,\F,P)\to\mathcal{X}$ and denote $X(B)$ by $X^B$ and $\mathscr{L}_B$ by 
${\rm dom }(X^B)$ for $B\in\pmb{B}_I(\Om,\F,P)$.
Let $h$ be a map over $\cup_{B\in\pmb{B}_I(\Om,\F,P)}\mathscr{L}_B$. 
For example, $h=\mathrm{id}\,($identity$),\, h=|\cdot|\,($absolute value or some norm$),\,E\,($expectation$)$ or $\,[\cdot]\,($some equivalence class$)$.
Let $W\in\pmb{B}_I(\Om,\F,P)$.

\begin{enumerate}
\i \textbf{Identification in the wide sense}\noindent

We say $h(\alpha)$ is identified  (in the wide sense) for $\alpha\in\mathscr{L}_W$
from $X^W(\alpha)=(X^W_\lambda(\alpha))_{\lambda\in\Lambda}$
if there exists a map $\mathcal{T}$ over ${\rm Im}(X^{W})$ such that 
$\mathcal{T}\circ X^{W}=h$ on $\mathscr{L}_W$.
In this case, we say $h(\alpha)$ is identified for $\alpha\in\mathscr{L}_W$ from 
$X^W(\alpha)$ by $\mathcal{T}$.\m

\i \textbf{Constructive identification in $\mac{L}$}\noindent

We say $h(\alpha)$ is identified constructively in $\mac{L}$ (or $\mac{L}$-constructively identified) for $\alpha\in\mathscr{L}_W$ from $X^W(\alpha)$
by the formula $h(\alpha)=s(X^{W}(\alpha))$ (or some abbreviation or notational deformation of this formula)
if $h(\alpha)$ is identified for $\alpha\in\mathscr{L}_W$ from $X^W(\alpha)$ by the map which is $\mac{L}$-constructive with the term $s(x)$.\m

\i \textbf{Constructive identification in $\mac{L}$ and $W$}\noindent

We say $h(\alpha)$ is identified constructively in $\mac{L}$ and (or with) $W$ (or $\mac{L}^W$-constructively identified) for $\alpha\in\mathscr{L}_W$ from $X^W(\alpha)$
by the formula $h(\alpha)=s(X^{W}(\alpha))$ (or some abbreviation or notational deformation of this formula)
if $h(\alpha)$ is identified for $\alpha\in\mathscr{L}_W$ from $X^W(\alpha)$ 
by the map which is $\mac{L}^W$-constructive with $s(x)$.\m

\i \textbf{B-independent identification}\noindent

We say $h(\alpha)$ is B-independently identified (identified in no need of (or without use of) information of $W$) for $\alpha\in\mathscr{L}_W$ from 
$X^W(\alpha)$
if there exists a map $\mathcal{T}$ such that for any $B\in\pmb{B}_I(\Om,\F,P)$ we have
\begin{equation}
{\rm Im}(X^{B})\subset{\rm dom}(\mathcal{T})\subset L^0(\Om)^{\Lambda}  \,\tx{ and }\,
\mathcal{T}(X^{B}(\beta))=h(\beta)  \,\tx{ for any }\beta\in\mathscr{L}_{B}.
\label{eq: 999}
\end{equation}

In this case, we say $h(\alpha)$ is B-independently identified for $\alpha\in\mathscr{L}_W$ from 
$X^W(\alpha)$ by $\mathcal{T}$.\m

\i \textbf{B-dependent identification}\noindent

We say $h(\alpha)$ is B-dependently identified (identified in need of (or with use of) information of $W$) for $\alpha\in\mathscr{L}_W$ from 
$X^W(\alpha)$
if $h(\alpha)$ is identified for $\alpha\in\mathscr{L}_W$ from $X^W(\alpha)$
and not B-independently identified for $\alpha\in\mathscr{L}_W$ from 
$X^W(\alpha)$.
\end{enumerate}
\label{def: 6.1}
\end{dfn}

\remi{1} The definitions of semantically constructive identifications are done in the same way as 
those of constructive identifications in the cases of 2 and 3.\\
\remi{2} In each of the cases above,
we say $\alpha\in\mathscr{L}_W$ is identified from $X^W(\alpha)$
if $\alpha$ is identified for $\alpha\in\mathscr{L}_{W}$ from $X^W(\alpha)$.\\
\remi{3} In the case of 4, if the map $\mac{T}$ is constructive in an assigned language $\mac{K}$ with the term $s(x)$,
we say $h(\alpha)$ is $\mac{K}$-constructively and B-independently identified for $\alpha\in\mathscr{L}_W$ from $X^W(\alpha)$ by the formula $h(\alpha)=s(X^{W}(\alpha))$ (or some abbreviation or notational deformation of this formula).\\
\remi{4} When $\Lambda\subset\mathbb{N}$
and $X^W_n(\alpha)$ is the $n$-th SFC, defined later, of $\alpha(t)$ for each $\alpha\in\mathscr{L}_W$,
the above definitions from 1 to 3 become those of
identifications of $h(\alpha)$ from SFCs of $\alpha\in\mathscr{L}_W$.\\
\remi{5} The identification of a stochastic differential (pair of random functions) can be defined similar to that of a random function.\\
\remi{6} The definition of B-independent identification is appropriate as a definition of identification in no need of the information that the underlying Brownian motion is $\s{W}{t}{I}$ as mentioned in the introduction by the following reason: 
the meaning of the sentence ''The derivation map $\mathcal{T}$ in the definition 1 is given
in no need of the information that the underlying Brownian motion is $\s{W}{t}{I}$'' is that 
we can give $\mathcal{T}$ before giving the underlying Brownian motion $W$ on $(\Om,\F,P)$.
So, such a map $\mathcal{T}$ must satisfy \eqref{eq: 999} for any Brownian motion $B$ on $(\Om,\F,P)$.\\
\remi{7} 
On the contrary,  B-dependent identification can be said to be identification in need of
the information that the underlying Brownian motion is $\s{W}{t}{I}$.\\
\remi{8} 
In light of Remarks 4 and 5,
B-independent identification can be said to be identification on the probability space $(\Om,\F,P)$
and B-dependent identification can be said to be proper identification on the probability space $(\Om,\F,P,\s{W}{t}{I})$ with $\s{W}{t}{I}$.\\
\remi{9} When $\mac{K}$ is an assigned first-order language without the function symbol which is assigned to the Brownian motion $W$
nor any other symbols which are assigned to elements depending on $W$,
$\mac{K}$-semantically constructive identification implies B-independent identification,
since a $\mac{K}$-semantically constructive map can be given before giving the underlying Brownian motion $W$ on $(\Om,\F,P)$.\\
\remi{10}
$h(\alpha)$ is B-independently identified for $\alpha\in\mathscr{L}_W$ from 
$X^W(\alpha)$ if and only if $X^{B^{(1)}}(\beta)=X^{B^{(2)}}(\gamma)$ implies $h(\beta)=h(\gamma)$
for any $B^{(1)},B^{(2)}\in\pmb{B}_I(\Om,\F,P)$ and $\beta\in\mathscr{L}_{B^{(1)}},\gamma\in\mathscr{L}_{B^{(2)}}$.\\


We now recall the setting given in Subsection 2.1.
From now on, let $(e_n)_{n\in\mab{N}}$ be a CONS of $L^2([0,L];\C)$, 
$a\in L^0([0,L]\t\Om\,;\C)$ and $b\in L^0(\Om; L^2([0,L];\C))$.

\begin{exm}
Let $(\vp_m)_{m\in\mab{N}}$ be a CONS of $L^2([0,L]\,;\,\mab{C})$
and
$f$ a non-zero random function in $L^0([0,L]\times\Omega)$.
Suppose $\overline{e_n}f$ is $\vp$-integrable with respect to $B$ for every $n\in\Lambda$, where $\emptyset\subsetneq\Lambda\subset\mathbb{N}$. 
Then, $a\in\{f,-f\}$ is not B-independently identified from its 
SFC-$\rm{O}_{\vp}$'s $((e_n,a\,d_{\vp}B))_{n\in\Lambda}$,
defined later,
since $B$ is symmetrically distributed.
Therefore, 
$a\in\{f,-f\}$ is not identified constructively (in both syntactic and semantic senses) in $\mac{K}$ in Remark 9 from 
$(\sfc{e_n}{a\,d_{\vp}B})_{n\in\Lambda}$.
In other words,
the sign of any random function cannot be identified by any map which can be defined before giving the underlying Brownian motion $B$ on $(\Om,\F,P)$.
\label{ex: 4.5}
\end{exm}

\begin{exm}
Let $\L_0$ be the assigned language defined in Example \r{ex:1}.
If $h(\alpha)$ in Definition \ref{def: 6.1} is identified constructively in $\mac{L}_0$ for $\alpha\in\mathscr{L}_B$ from $X^B(\alpha)$
by the formula $h(\alpha)=s(X^{B}(\alpha))$,
$h(\alpha)$ is B-independently identified for $\alpha\in\mathscr{L}_B$ from $X^B(\alpha)$ by the same formula,
since $\L_0$ satisfies the assumption in Remark 9.
\label{ex:4.6}
\end{exm}

\begin{exm}
Let $\L_0$ be the assigned language defined in Example \r{ex:1}.
Then
in Theorem \r{thm: 0001},
$a(t)$ is $\L_0$-constructively
(thus, and B-independently)
identified from $(\hat{a}_n)_{n\in\mab{Z}}$ by \eqref{al:0.00100}.
\end{exm}


\section{Lemmas}

In this section, we state primary lemmas necessary to obtain main theorems.
First, the next lemma follows from Doob's $L^2$-inequality.

\begin{lem}\bt
Let $(f_n)_{n\in\mab{N}}$ be a sequence of $L^2([0,L];\mab{C})$ which converges to 0.
Then, the following holds:
$$
\lim_{n\to\infty}\sup_{s\in[0,L]}|B_s[f_n]|=0 \q\tx{in}\,\, L^2(\Om).
$$
\label{lem:4}
\end{lem}

\begin{prf}\bt
Fix $n\in\mab{N}$. Since $(B_s[f_n])_{s\in[0,L]}$ is a martingale with respect to the filtration generated by $(B_s)_{s\in[0,L]}$,
then by Doob's $L^2$-inequality we have 
\begin{align*}
E\Bigl(\,\Bigl(\sup_{\,s\in[0,L]}|B_s[f_n]|\,\Bigl)^2\,\Bigl)
\le&\,4\sup_{s\in[0,L]} E|B_s[f_n]|^2.
\end{align*}
Besides, by the isometry of the \wi integral we have 
\begin{align}
\sup_{s\in[0,L]} E|B_s[f_n]|^2
=&\,\int_0^L|f_n|^2\,d\lambda
\label{al:4.003}
\end{align}
and \eqref{al:4.003} converges to $0$ as $n\to\infty$.
Therefore we have
\mbox{
$
\tlim_{n\to\infty}\!E\Bigl(\,\Bigl(\tsup_{\,s\in[0,L]}|B_s[f_n]|\,\Bigl)^2\,\Bigl)=0.
$
}
\end{prf}

The next lemma follows from Lemma \r{lem:4} and Proposition \r{ap:24}.

\begin{lem} 
Assume $a(t)$ is a \nc finite variation process on $[0,L]$.
Let $(f_n)_{n\in\mab{N}}$ be a sequence of $L^2([0,L];\mab{C})$ such that $\tlim_{n\to\infty}f_n=0 \, \tx{ in } L^2([0,L];\mab{C})$.
Then, the following hold:
$$
\lim_{n\to\infty}\int_{(0,L)}B_s[f_n]\,d\,\widetilde{a}(s)=0 \q\tx{in probability},
$$
where $\widetilde{a}(t)$ is the version of $a(t)$ given by Proposition \r{ap:24}.
\label{lem:4.1}
\end{lem}
\begin{prf}
For each $n\in\mab{N}$ we have
\begin{align}
\Bigl|\,\int_{(0,L)}B_s[f_n]\,d\,\widetilde{a}(s)\,\Bigl|
\le&\,||\,\widetilde{a}\,||\sup_{s\in[0,L]}|B_s[f_n]|.
\label{al:4.004}
\end{align}
Here, $||\,\widetilde{a}\,||$ is a random variable by Proposition \r{ap:24} and $||\,\widetilde{a}\,||<\infty$ a.s.
Then, by Lemma \r{lem:4}, the right hand side of \eqref{al:4.004} converges to 0 in probability, which completes the proof of Lemma \r{lem:4.1}.
\end{prf}


\section{Identification of random functions from \sfcoS}

Hereafter, let $\L_0$ be the assigned first-order language defined in Example \r{ex:1}.
Hereafter, we also mean by a term (resp. formula) of $\mac{K}$ some abbreviation or notational deformation of a term (resp. formula) of $\mac{K}$ for an assigned first-order language.
In this section,
we give the results about identification of random functions from \sfcos of a stochastic differential by $a(t)$ and $b(t)$.

\subsection{Identification of \nc finite variation processes}
\label{sub: 6.1}

In this subsection, we give the main results about identification from \sfcos of a stochastic differential whose diffusion coefficient is any \nc finite variation process.

Assume $a\in L^0([0,L]\t\Om)$ is any \nc finite variation process.
First, we give a necessary and sufficient condition for a random function or some quantity to be identified from SFC-Os.

\begin{prp}[Parseval-type transformation]
Assume $e_n(t)$ has compact support for each $n\in\mab{N}$.
Then, the following hold:
\begin{enumerate}
\i[(a)] $\tlim_{n\to\infty}\int_0^Lf_na\,d_{\rm{u}}B=0$ in probability for any sequence $(f_n)_{n\in\mab{N}}$ of $L^2([0,L];\mab{C})$ which converges to 0 such that $\a n\in\mab{N} \,\, f_n$ has compact support.
\end{enumerate}
In particular,
\begin{itemize}
\i[(b)] \inum
\begin{align}
\!\!\!\!\!\!\!\!\!\!\!\!\!\!\mathcal{P}((e_n,a\,d_{\rm u}B))_{n\in\mab{N}}(t):
&=\underset{N\to\infty}{\rm{l.i.p.}}\sum_{n=1}^{N}\int_0^te_n(s)\,ds\,(e_n,a\,d_{\rm u}B)\nonumber\\
&=\displaystyle\int_0^ta\,d_{\rm u}B,\,\,\a t\in[0,L],
\label{al:4.11152500}
\end{align}
\i[(c)] \inum
\begin{align}
\!\!\!\!\!\!\!\!\!\!\!\!\!\!\mathcal{P}((e_n,d_{\rm u}Y))_{n\in\mab{N}}(t):
&=\underset{N\to\infty}{\rm{l.i.p.}}\sum_{n=1}^{N}\int_0^te_n(s)\,ds\,(e_n,d_{\rm u}Y)\nonumber\\
&=Y_t,\,\,\a t\in[0,L],
\label{al:4.111525}
\end{align}
\end{itemize}
where $d_{\rm u}Y$ denotes the stochastic differential $d_{\rm u}Y_t=a(t)\,d_{\rm u}B_t+b(t)\,dt$
and $\mathcal{P}$ is $\mac{L}_0$-constructive.
\label{pr:4.112}
\end{prp}
\rem  We can extend $\mac{P}$ for $((e_n,d_{\rm u}Y))_{n\in\Lambda}$, where $\Lambda$ is a cofinite subset of $\mab{N}$,
by regarding $(e_n,d_{\rm u}Y)=0$ if $n\notin\Lambda$.

\begin{prf} 
First of all, \sfcO s $\,\sfc{e_n}{a\dub}$ and $\int_0^Lf_na\,d_{\rm u}B$ are well-defined because $\overline{e_n}a$ and $f_na$ is u-integrable for $L^2([0,L]\,;\mab{R})$ by Lemma \r{lem:0.7}.
By this lemma,
the \og integral of $f_na$ is given by 
\begin{align}
\int_0^L
	f_na
  \,d_{\rm u}B
=\tilde{a}(L)B_L[f_n]-\int_{(0,L)}B_s[f_n]\,d\tilde{a}(s)\q\tx{almost surely,}
\label{al:1.5578}
\end{align}
where $\tilde{a}(t)$ is the version of $a(t)$ as in Proposition \r{ap:24}. 
The convergence as $n\to\infty$ of the first term of the right hand side of \eqref{al:1.5578} is that of \wi expansion to 0, and the second term converges to 0 in probability by Lemma \r{lem:4.1}. Therefore (a) holds.
(b) is obtained by taking $\ic{[0,t]}-\tsum_{i=1}^{n}\overline{\ip{e_i}{\mathsf{1}_{[0,t]}}{}e_i}$ for $f_n$ for each $t\in[0,L]$. (b)$\Rightarrow$(c) is obvious.
\end{prf} 

\begin{cor}
Assume $e_n(t)$ has compact support for each $n\in\mab{N}$.
Then, for a subset $\mathscr{S}$ of the family of all \nc finite variation processes on $[0,L]$ and a dense subset $S$ of $[0,L]$ and a map $h$ over $\mathscr{S}$,
the following are equivalent:
\begin{enumerate}
\i[$(i)$] $h(a)$ is identified for $a\in\mathscr{S}$ from $((e_n,d_{\rm{u}}Y))_{n\in\mab{N}}$.
\i[$(ii)$] $h(a)$ is identified for $a\in\mathscr{S}$ from $(Y_t)_{t\in S}$.
\end{enumerate}
Here $d_{\rm u}Y$ denotes the stochastic differential $d_{\rm u}Y_t=a(t)\,d_{\rm u}B_t+b(t)\,dt$.
\label{cor:4.112}
\end{cor}

\begin{prf} 
The implication $(ii)\Rightarrow(i)$ is obtained from (c) in Proposition \ref{pr:4.112}.
On the other hand, $(i)\Rightarrow(ii)$ is justified because each $e_n(t)$ is approximated in $L^2[0,L]$ by step functions of the form $\tsum_{j=1}^{r}\alpha_j\ic{[0,t_j]}$, $r\in\mab{N}$, $\alpha_1,\ldots,\alpha_r\in\mab{C}$ and $t_1,\ldots,t_r\in S$.
\end{prf} 

\begin{lem}
Let $\Lambda$ be a cofinite subset of $\mab{N}$.
Suppose $e_n(t)$ has compact support for each $n\in\Lambda$.
Then,
the map 
which associates $\mathcal{P}((e_n,d_{\rm u}Y))_{n\in\Lambda}(t)$ in Proposition \ref{pr:4.112} with its left-continuous modification $\widetilde{\mathcal{P}}((e_n,d_{\rm u}Y))_{n\in\Lambda}(t)$ is given by
\[
\widetilde{\mathcal{P}}((e_n,d_{\rm u}Y))_{n\in\Lambda}(t)
=\begin{cases}
\tlim_{\substack{s\nearrow t\\s\in\Q}}\mathcal{P}((e_n,d_{\rm u}Y))_{n\in\Lambda}(s)&,t>0\\
0&,t=0.
\end{cases}
\]
Then, the map is $\mac{L}_0$-constructive.
\label{lem:6.3}
\end{lem}
\begin{prf}
Consider the case $\Lambda=\mab{N}$ for simplicity.
By integration by parts we have
\begin{align}
\mathcal{P}((e_n,d_{\rm u}Y))_{n\in\mab{N}}(t)
&=\int_0^ta\,d_{\rm u}B+\int_0^tb\,d\lambda\nonumber\\
&=\tilde{a}(t)B_t-\int_{(0,t)}B\,d\tilde{a}+\int_0^tb\,d\lambda
\q \tx{almost surely.}
\label{al:1.5654}
\end{align}
Note that \eqref{al:1.5654} holds for each fixed $t\in[0,L]$.
Then, we have
$$
\Bigl(\,\mathcal{P}((e_n,d_{\rm u}Y))_{n\in\mab{N}}(t)=\tilde{a}(t)B_t-\int_{(0,t)}B\,d\tilde{a}+\int_0^tb\,d\lambda\q\a t\in\Q\cap[0,L]\,\Bigl)\,\, \tx{a.s.}
$$
and the process $\bigl(\,\tilde{a}(t)B_t-\int_{(0,t)}B\,d\tilde{a}+\int_0^tb\,d\lambda\,\bigl)_{t\in[0,L]}$ is left-continuous almost surely. 
Then, we can define a process $(\widetilde{X}_t)_{t\in[0,L]}$ pathwisely by
\[
\widetilde{X}_t:
=\begin{cases}
\tlim_{\substack{s\nearrow t\\s\in\Q}}\mathcal{P}((e_n,d_{\rm u}Y))_{n\in\mab{N}}(s)&,t>0\\
0&,t=0
\end{cases}
\]
which is equal to $\bigl(\,\tilde{a}(t)B_t-\int_{(0,t)}B\,d\tilde{a}+\int_0^tb\,d\lambda\,\bigl)_{t\in[0,L]}$ pathwisely almost surely.
Thus, $\s{\widetilde{X}}{t}{[0,L]}$ is a left-continuous realization of $\mathcal{P}((e_n,d_{\rm u}Y))_{n\in\mab{N}}(t)$ 
and the formula $s\in\Q$ is equivalent to the formula $g(s)=1$ of $\L_0$
with the Dirichlet function $g$, which is $\mac{L}_0$-constructive (Example \r{ex:1}(d)), which implies the assertion of this lemma.
\end{prf}

Hereafter, we give the main theorems.
\vspace{8pt}

\noindent
$\bullet\,$ {\bf Identification from SFCs of $\boldsymbol{a(t)}$}
\vspace{8pt}

At first, we give several reconstruction formulas from SFCs $((e_n,a\,d_{\rm u}B))_{n\in\mab{N}}$.

\begin{thm}[Constructive identification]
Let $\Lambda$ be a cofinite subset of $\mab{N}$.
Suppose $e_n(t)$ has compact support for each $n\in\Lambda$ and
$a\in L^0([0,L]\t\Om)$ with a \nc finite variation process as a version.
Then, the following hold:
\begin{align}
\,\,\displaystyle\Biggl(\,|a|(t)=\lims{s\searrow t}\frac{\,\widetilde{\mathcal{P}}((e_n,a\,d_{\rm u}B))_{n\in\Lambda}(s)-\widetilde{\mathcal{P}}((e_n,a\,d_{\rm u}B))_{n\in\Lambda}(t)}{\sqrt{2(s-t)\log\log\frac{1}{s-t}}} \,\,\, \tx{a.s.} \,\Biggl)\,\, \tx{a.a.} \,\,t\in[0,L],
\label{al:4.11154}
\end{align}
where $\widetilde{\mathcal{P}}$ is defined
in Lemma \ref{lem:6.3}.
Therefore,
$|a|$ is $\mac{L}_0$-constructively (thus, and B-independently) identified from $((e_n,a\,d_{\rm u}B))_{n\in\Lambda}$ by \eqref{al:4.11154},
where $\L_0$ is the assigned first-order language defined in Example \r{ex:1}.
\label{thm:6.20}
\end{thm}

\begin{prf}
Consider the case of $\Lambda=\mab{N}$ first.
By the definition of $\widetilde{X}_t=\widetilde{\mathcal{P}}((e_n,a\,d_{\rm u}B))_{n\in\mab{N}}(t)$ in the case of $b=0$ in Lemma \ref{lem:6.3}, we have
\begin{equation}
\Bigl(\,\widetilde{X}_t=\tilde{a}(t)B_t-\int_{(0,t)}B\,d\tilde{a}\q\a t\in[0,L]\,\Bigl)\,\, \tx{a.s.}
\label{eq:4.02}
\end{equation}
On the other hand, because $\tilde{a}(t)$ is of bounded variation,
we have the following by \leB's theorem:
\begin{equation}
\tx{( } \tilde{a}(\cdot,\om) \tx{ and } \tilde{a}_{\rm tv}(\cdot,\om) \tx{ is differentiable in } t\,\tx{ a.a.} \,\, t\in[0,L] \tx{ ) \, a.a.} \,\, \om\in\Om.
\label{eq:4.1}
\end{equation}
By Proposition \r{ap:24.5}, \eqref{eq:4.1} is rephrased as
\begin{equation}
\tx{( } \tilde{a}(\cdot,\om) \tx{ and } \tilde{a}_{\rm tv}(\cdot,\om) \tx{ is differentiable in } t \,\tx{ a.a.} \,\, \om\in\Om \tx{ ) \, a.a.} \,\, t\in[0,L].
\label{eq:4.03}
\end{equation}
Now, by \eqref{eq:4.02} and \eqref{eq:4.03} for almost all $t\in[0,L]$, taking account of laws of the iterated logarithm  
$\overline{\tlim_{s\searrow t}}\frac{B_s-B_t}{\sqrt{2(s-t)\log\log\frac{1}{s-t}}}=1 \tx{ and }\,\,\underset{s\searrow t}{\underline{\tlim}} \frac{B_s-B_t}{\sqrt{2(s-t)\log\log\frac{1}{s-t}}}=-1 \tx{ a.s.} $,
we have 
\begin{align*}
c(t):
=&\overline{\lim_{s\searrow t}}\frac{\widetilde{X}_s-\widetilde{X}_t}{\sqrt{2(s-t)\log\log\frac{1}{s-t}}}
=\overline{\lim_{s\searrow t}}\,\frac{\tilde{a}(s)B_s-\tilde{a}(t)B_t-\int_{[t,s)}B\,d\tilde{a}}{\sqrt{2(s-t)\log\log\frac{1}{s-t}}}\\
=&\overline{\tlim_{s\searrow t}}\Biggl(\frac{B_{s}-B_t}{\sqrt{2(s-t)\log\log\frac{1}{s-t}}}\,\tilde{a}(t)+\frac{\tilde{a}(s)-\tilde{a}(t)}{s-t}\,\frac{\sqrt{s-t}B_{s}}{\sqrt{2\log\log\frac{1}{s-t}}}-\frac{\int_{[t,s)}B\,d\tilde{a}}{\sqrt{2(s-t)\log\log\frac{1}{s-t}}}\Biggl)
\end{align*}
almost surely. Here
\begin{align*}
\overline{\lim_{s\searrow t}}\frac{B_{s}-B_t}{\sqrt{2(s-t)\log\log\frac{1}{s-t}}}\,\tilde{a}(t)=|\tilde{a}(t)|\q \tx{and}\q
\lim_{s\searrow t}\,\frac{\tilde{a}(s)-\tilde{a}(t)}{s-t}\,\frac{\sqrt{s-t}B_{s}}{\sqrt{2\log\log\frac{1}{s-t}}}=\tilde{a}'(t)\cdot 0=0
\end{align*}
almost surely. Moreover when $s$ $(>t)$ is close to $t$, we have
\begin{align*}
\Biggl|\,\frac{\int_{[t,s)}B\,d\tilde{a}}{\sqrt{2(s-t)\log\log\frac{1}{s-t}}}\,\Biggl|
&\le\frac{(\tilde{a}_{\rm tv}(s)-\tilde{a}_{\rm tv}(t))\max_{u\in[t,t+1]}|B_u|}{\sqrt{2(s-t)\log\log\frac{1}{s-t}}}\nonumber\\
&=\frac{\tilde{a}_{\rm tv}(s)-\tilde{a}_{\rm tv}(t)}{s-t}\frac{\sqrt{s-t}\max_{u\in[t,t+1]}|B_u|}{\sqrt{2\log\log\frac{1}{s-t}}},
\end{align*}
and the right hand side converges to $\tilde{a}_{\rm tv}'(t)\cdot 0=0$ as $s\searrow t$ almost surely.
Then, we have 
\begin{align}
(\,c(t)=|\tilde{a}(t)| \,\tx{ a.s. }\,)\, \tx{ a.a. }\,t\in[0,L].
\label{al:16.12}
\end{align}
Here,
$c\in L^0([0,L]\t\Om)$ by the following reason:
Putting 
$f(t,s,\om)=\frac{\widetilde{X}_s(\om)-\widetilde{X}_t(\om)}{\sqrt{2(s-t)\log\log\frac{1}{s-t}}}$ and 
$g_n(t,\om)=\tsup_{t<s<t+\frac{1}{n}}f(t,s,\om),\,n\in\N$,
then
$$
(\,\a t\in[0,L]\,\,f(t,s)\tx{ is left-continuous at }s\in(t,L]\,)\tx{ a.s.}
$$
So, there exists $\tilde{\Om}\in\mac{F}$ such that $P(\tilde{\Om})=1$ and
\begin{align*}
\!\{(t,\om)\in[0,L]\t\tilde{\Om}\,|\,g_n(t,\om)\!>r\,\}
&=\{\,(t,\om)\in[0,L]\t\tilde{\Om}\,|\,\e s\in(t,t+\frac{1}{n})\,\,f(t,s,\om)>r\,\}\\
&=\{\,(t,\om)\in[0,L]\t\tilde{\Om}\,|\,\e s\in(t,t+\frac{1}{n})\cap\Q\,\,f(t,s,\om)>r\,\}\\
&=\!\cu_{s\in\Q}\Bigl((s-\frac{1}{n},s)\!\t\!\tilde{\Om}\,\,\cap\,\, f(\cdot,s,\cdot)^{-1}((r,\infty])\Bigl)\,\in\mac{L}([0,L])\!\ot\!\mac{F}
\end{align*}
for any $r\in\R$.
Then,
$c=\tlim_{n\to\infty}g_n\in L^0([0,L]\t\Om)$.
Therefore,
we see from (\ref{al:16.12}) that $c=|a| \,\,\lambda\ot\P\tx{-a.e.}$ by Fubini's theorem.
Finally, the formula \eqref{al:4.11154}
in the case of $\Lambda\subsetneq\mab{N}$ is obtained because for  all $t\in[0,L]$ and $n\in\Lambda^c$ almost surely we have by \sci
\begin{align*}
\Biggl|\,\frac{\overline{\ip{e_n}{\ic{[0,s]}}{}}\sfc{e_n}{a\dub}^{\widetilde{}}-\overline{\ip{e_n}{\ic{[0,t]}}{}}\sfc{e_n}{a\dub}^{\widetilde{}}}{\sqrt{2(s-t)\log\log\frac{1}{s-t}}}\,\Biggl|\,
\le\,&\frac{|\sfc{e_n}{a\dub}^{\widetilde{}}\,|}{\sqrt{2\log\log\frac{1}{s-t}}},
\end{align*}
and the left hand side converges to 0 as $s\searrow t$ since the right hand side converges to 0.
Here $\sfc{e_n}{a\dub}^{\widetilde{}}$ denotes $\tilde{a}(L)B_L[\overline{e_n}]-\int_{(0,L)}B_s[\overline{e_n}]\,d\tilde{a}(s)$, noting that 
$\sfc{e_n}{a\dub}$ is not necessarily defined for $n\in\Lambda^c$.
\end{prf}

\begin{thm}[Constructive identification with $B$]
Let $\Lambda$ be a cofinite subset of $\mab{N}$.
Suppose $e_n(t)$ has compact support for each $n\in\Lambda$ and
$a\in L^0([0,L]\t\Om)$ with a \nc finite variation process as a version.
Then, the following hold:
\begin{align}
\,\,\displaystyle\Biggl(\,\, a(t)
=\lim_{k\to\infty}\Biggl(\,\,\lims{s\searrow t}\,\frac{\widetilde{\mathcal{P}}((e_n,a\,d_{\rm u}B))_{n\in\Lambda}(s,t)+k(B_s-B_t)}{\sqrt{2(s-t)\log\log\frac{1}{s-t}}}-k\,\,\Biggl) \q \tx{a.s.}
\,\,\Biggl)\,\,\, \tx{a.a.} \,\,t\in[0,L],
\label{al:4.11155}
\end{align}
where $\widetilde{\mathcal{P}}$ is defined
in Lemma \ref{lem:6.3}
and $\widetilde{\mathcal{P}}((e_n,a\,d_{\rm u}B))_{n\in\Lambda}(s,t)$ denotes
$\widetilde{\mathcal{P}}((e_n,a\,d_{\rm u}B))_{n\in\Lambda}(s)-\widetilde{\mathcal{P}}((e_n,a\,d_{\rm u}B))_{n\in\Lambda}(t)$.
Therefore, $a(t)$ is  identified constructively in $\mac{L}_0$ and $B$ from $((e_n,a\,d_{\rm u}B))_{n\in\Lambda}$ by (\r{al:4.11155}),
where $\L_0$ is the assigned first-order language defined in Example \r{ex:1}.
\label{pr:4.11466}
\end{thm}

\begin{prf}
Consider the case of $\Lambda=\mab{N}$.
Put $a_k:=a+k$ for each $k\in\mab{N}$. 
Notice $\int_t^s a\,d_{\rm u}B+k(B_s-B_t)=\int_t^sa_k\,d_{\rm u}B$ and
$\widetilde{\mathcal{P}}((e_n,a\,d_{\rm u}B))_{n\in\mab{N}}(s,t)+k(B_s-B_t)$
equals
$\widetilde{\mathcal{P}}((e_n,a_k\,d_{\rm u}B))_{n\in\mab{N}}(s,t)$ as left-continuous processes.
By (\ref{al:4.11154}), we have to check
\begin{align}
\bigl(\,\,a(t,\om)=\lim_{k\to\infty}\bigl(\,|a_k(t,\om)|-k\,\bigl)\q\tx{a.a.}\,\,\om\in\Om\,\,\bigl)\q\tx{a.a.}\,\, t\in[0,L].
\label{al:4.0451}
\end{align}
For almost all $\om\in\Om$, $a(\cdot,\om)$ is bounded, so there exists $k(\om)\in\mab{N}$ such that  $\a k\ge k(\om)\,\,\a t\in[0,L]\,\,\,\tilde{a}_k(t,\om)\ge0$. Then, \eqref{al:4.0451} holds.
Besides,
the right hand side of \eqref{al:4.11155} is defined
by using the function symbol $B$ for the term $k(B_s-B_t)$.
Then, $a(t)$ is identified constructively in $\mac{L}_0$ and $B$ from $((e_n,a\,d_{\rm u}B))_{n\in\mab{N}}$.
Similarly, the proof
in the case of $\Lambda\subsetneq\mab{N}$ is done.
\end{prf}

\begin{cor}
Let $\L_0$ be the assigned first-order language defined in Example \r{ex:1},
$\Lambda$ a cofinite subset of $\mab{N}$ and $T=[0,L)$ or $\{0\}$.
Suppose $e_n(t)$ has compact support for each $n\in\Lambda$.
Let $\hat{a}(t)$ be a \nc finite variation process on $[0,L]$
which is right-continuous in $T$ and left-continuous in $[0,L]\diagdown T$ with probability 1.
Put $a=[\hat{a}]\in L^0([0,L]\t\Om)$.
Then, the following hold:
\begin{enumerate}
\i[(A)] $|\hat{a}|$ is $\mac{L}_0$-constructively (thus, and B-independently) identified almost surely from $(\sfc{e_n}{a\dub})_{n\in\Lambda}$
by \eqref{al:4.11154} and
\begin{equation*}
|\hat{a}|(t)=
\begin{cases}
\lims{n\to\infty}n\int_t^{t+\frac{1}{n}} |a|\,d\lambda	&,t\in T\\
\lims{n\to\infty}n\int_{t-\frac{1}{n}}^t |a|\,d\lambda	&,t\notin T.
\end{cases}
\end{equation*}
\i[(B)] $\hat{a}$ is  identified constructively in $\mac{L}_0$ and $B$ almost surely from $(\sfc{e_n}{a\dub})_{n\in\Lambda}$
by \eqref{al:4.11155} and
\begin{equation*}
\hat{a}(t)=
\begin{cases}
\lims{n\to\infty}n\int_t^{t+\frac{1}{n}} a\,d\lambda	&,t\in T\\
\lims{n\to\infty}n\int_{t-\frac{1}{n}}^t a\,d\lambda	&,t\notin T.
\end{cases}
\end{equation*}
\end{enumerate}
\label{thm: 1.5}
\end{cor}

\remi{1} Almost surely (omitted below), 
$|\hat{a}|(t)$ (so is $\hat{a}(t)$) is identified for every $t\in[0,L]$,
that is, identified as a function with the equivalence
'' $\hat{b}=\hat{c} \Leftrightarrow (\,\hat{b}(t)=\hat{c}(t)\,\,\a t\in[0,L]\,)$ a.s.''.
So,
$|a|(t)$ (so is $a(t)$) is identified for every $t\in[0,L]$ except differences of countable discontinuous points. 
But it does not mean almost surely, $|a|(t)$ (resp. $a(t)$) is identified for every $t\in [0,L]$, because two finite variation processes equal at almost every where $(t,\om)\in[0,L]\t\Om$ define same \sfcO.\\
\remi{2} 
The same assertion holds, even if the set $T$ in this corollary is $\{\,t\in[0,L]\,|\,{\rm eval}_{V}^{\sigma_t}$
$(\vp(x))_1\,\}\in\mac{L}([0,L])$
for some formula $\vp(x)$ of $\mac{L}_0$ with one variable $x$
such that ''${\rm syn}_{V}^{\sigma}(\vp(x))_2=1$ for any variable assignment $\sigma:\{x\}\to [0,L]$'' is provable,
and $\sigma_t:\{x\}\to \{t\}$.

\begin{prf} 
We assume $T=[0,L)$ for simplicity.
To prove (A),
it suffices to show $|\hat{a}|$ is identified from
the function $c(t)$ obtained in the proof of Theorem \ref{thm:6.20}.
First, we have
$(\,c(t)=|\hat{a}(t)| \,\tx{ a.a. }\,t\in[0,L]\,)\, \tx{ a.s.}$
Now, we set
$$
\tilde{c}(t):=\lims{n\to\infty}n\int_{t}^{t+\frac{1}{n}} c\,d\lambda,\,t\in[0,L).
$$
Then almost surely, we have
\begin{align*}
\tilde{c}(t)&=\lim_{n\to\infty}n\int_{t}^{t+\frac{1}{n}}  |\hat{a}|\,d\lambda\\
&=|\hat{a}(t)|, \,\,\a t\in[0,L)
\end{align*}
since $|\hat{a}|$ is right-continuous almost surely.
Similarly, $\tilde{c}(L)$ is obtained, or
$\tilde{c}(t)$ can be extend to $[0,L]$ left-continuously at $L$.
This completes the proof of (A).
(B) is proved by the same argument as in the proof of (A).
\end{prf}

\vspace{8pt}

\noindent
$\bullet\,$ {\bf Identification from SFCs of $\boldsymbol{d_{\rm u}Y=a(t)\,d_{\rm u}B_t+b(t)\,dt}$}
\vspace{8pt}

Next, we give the cores of the main theorems which give several reconstruction formulas from SFCs $((e_n,d_{\rm u}Y))_{n}$
as natural extensions to stochastic differentials
of Theorem \r{thm:6.20}, \r{pr:4.11466} and Corollary \r{thm: 1.5}.

\begin{thm}[Constructive identification]
Let $\Lambda$ be a cofinite subset of $\mab{N}$.
Suppose $e_n(t)$ has compact support for each $n\in\Lambda$  and
$a\in L^0([0,L]\t\Om)$ with a \nc finite variation process as a version.
Then, letting $d_{\rm u}Y=a(t)\,d_{\rm u}B_t+b(t)\,dt$, $b\in L^0(\Omega\,;L^2([0,L]\,;\C))$ the following hold:
\begin{align}
\,\,\displaystyle\Biggl(\,\,|a|(t)=\lims{s\searrow t}\frac{\,\widetilde{\mathcal{P}}((e_n,d_{\rm u}Y))_{n\in\Lambda}(s)-\widetilde{\mathcal{P}}((e_n,d_{\rm u}Y))_{n\in\Lambda}(t)}{\sqrt{2(s-t)\log\log\frac{1}{s-t}}} \q \tx{a.s.} \,\,\Biggl)\q \tx{a.a.} \,\,t\in[0,L],
\label{al:4.111545}
\end{align}
where $\widetilde{\mathcal{P}}$ is defined
in Lemma \ref{lem:6.3}.
Therefore, $|a|$ is $\mac{L}_0$-constructively (thus, and B-independently) identified from $((e_n,d_{\rm{u}}Y))_{n\in\Lambda}$ by \eqref{al:4.111545},
where $\L_0$ is the assigned first-order language defined in Example \r{ex:1}.
\label{thm:6.220}
\end{thm}

\begin{prf}
Consider the case $\Lambda=\N$.
By the definition of $\widetilde{X}_t=\widetilde{\mathcal{P}}((e_n,d_{\rm u}Y))_{n\in\mab{N}}(t)$ in Lemma \ref{lem:6.3}, we have
\begin{equation*}
\Bigl(\,\widetilde{X}_t=\tilde{a}(t)B_t-\int_{(0,t)}B\,d\tilde{a}+\int_0^tb\,d\lambda\q\a t\in[0,L]\,\Bigl)\,\, \tx{a.s.}
\end{equation*}
Thus, as in the proof of Theorem \ref{thm:6.20} we gain
$\Bigl(\,c(t):=\overline{\tlim_{s\searrow t}}\frac{\widetilde{X}_s-\widetilde{X}_t}{\sqrt{2(s-t)\log\log\frac{1}{s-t}}}=|a|(t) \q\tx{a.s.}\,\Bigl)$
$\,\,\tx{a.a.}\,\, t\in[0,L]
$
because
\begin{align}
\Biggl|\,\frac{\int_t^sb\,d\lambda}{\sqrt{2(s-t)\log\log\frac{1}{s-t}}}\,\Biggl|
\le\frac{(\int_t^s|b|^2d\lambda)^{\frac{1}{2}}\sqrt{s-t}}{\sqrt{2(s-t)\log\log\frac{1}{s-t}}}
&\le\frac{|b|_{L^2([0,L]\,;\mab{C})}}{\sqrt{2\log\log\frac{1}{s-t}}},
\label{al:20.34}
\end{align}
which shows the left hand side of \eqref{al:20.34} converges to 0 as $s\searrow t$.\\
Here, $c\in L^0([0,L]\t\Om)$ and $c=|a| \,\,\lambda\ot\P\tx{-a.e.}$ as in the proof of Theorem \ref{thm:6.20}.
\end{prf}

\rem \eqref{al:4.111545} also gives the law of iterated logarithm for the left-continuous modification of the process $Y_t$.

\begin{thm}[Constructive identification with $B$]
Let $\Lambda$ be a cofinite subset of $\mab{N}$.
Suppose $e_n(t)$ has compact support for each $n\in\Lambda$ and
$a\in L^0([0,L]\t\Om)$ with a \nc finite variation process as a version.
Then, letting $d_{\rm u}Y=a(t)\,d_{\rm u}B_t+b(t)\,dt$, $b\in L^0(\Omega\,;L^2([0,L]\,;\C))$ the following hold:
\begin{align}
\,\,\displaystyle\Biggl(\,\,a(t)=\lim_{k\to\infty}\Biggl(\,\lims{s\searrow t}\frac{\,\widetilde{\mathcal{P}}((e_n,d_{\rm u}Y))_{n\in\Lambda}(s,t)+k(B_s-B_t)}{\sqrt{2(s-t)\log\log\frac{1}{s-t}}}-k\,\Biggl) \q \tx{a.s.} \,\,\Biggl)\q \tx{a.a.} \,\,t\in[0,L],
\label{al:4.111555}
\end{align}
where $\widetilde{\mathcal{P}}$ is defined
in Lemma \ref{lem:6.3}
and $\widetilde{\mathcal{P}}((e_n,d_{\rm u}Y))_{n\in\Lambda}(s,t)$ denotes
$\widetilde{\mathcal{P}}((e_n,d_{\rm u}Y))_{n\in\Lambda}(s)-\widetilde{\mathcal{P}}((e_n,d_{\rm u}Y))_{n\in\Lambda}(t)$.
Therefore, $a(t)$ is identified  constructively in $\mac{L}_0$ and $B$ from $((e_n,d_{\rm{u}}Y))_{n\in\Lambda}$ by \eqref{al:4.111555},
where $\L_0$ is the assigned first-order language defined in Example \r{ex:1}.
\label{pr:4.11440}
\end{thm}


\begin{prf}
The proof is done in the same way of that of Theorem \ref{pr:4.11466} by using \eqref{al:4.111545}. 
\end{prf}


\begin{cor} 
Let $\L_0$ be the assigned first-order language defined in Example \r{ex:1},
$\Lambda$ a cofinite subset of $\mab{N}$ and $T=[0,L)$ or $\{0\}$.
Suppose $e_n(t)$ has compact support for each $n\in\Lambda$.
Let $\hat{a}(t)$ be a \nc finite variation process on $[0,L]$
which is right-continuous in $T$ and left-continuous in $[0,L]\diagdown T$ with probability 1. 
Put $a=[\hat{a}]\in L^0([0,L]\t\Om)$.
Then, letting $d_{\rm u}Y=a(t)\,d_{\rm u}B_t+b(t)\,dt$, $b\in L^0(\Omega\,;L^2([0,L]\,;\C))$, $b^{\Lambda}=\tsum_{n\in\Lambda}\ip{e_n}{b}{}e_n\in L^0(\Omega\,;L^2([0,L]))$,
the following hold:
\begin{enumerate}
\i[(A)] 
$|\hat{a}|$ is $\mac{L}_0$-constructively (thus, and B-independently) identified almost surely from $(\sfc{e_n}{d_{\rm u}Y})_{n\in\Lambda}$ by
\eqref{al:4.111545} and
\begin{equation*}
|\hat{a}|(t)=
\begin{cases}
\lims{n\to\infty}n\int_t^{t+\frac{1}{n}} |a|\,d\lambda	&,t\in T\\
\lims{n\to\infty}n\int_{t-\frac{1}{n}}^t |a|\,d\lambda	&,t\notin T.
\end{cases}
\end{equation*}
\i[(B)] 
$\hat{a}$ is identified constructively in $\mac{L}_0$ and $B$ almost surely from $(\sfc{e_n}{d_{\rm u}Y})_{n\in\Lambda}$ by
\eqref{al:4.111555} and
\begin{equation*}
\hat{a}(t)=
\begin{cases}
\lims{n\to\infty}n\int_t^{t+\frac{1}{n}} a\,d\lambda		&,t\in T\\
\lims{n\to\infty}n\int_{t-\frac{1}{n}}^t a\,d\lambda		&,t\notin T.
\end{cases}
\end{equation*}
\i[(C)] $b^{\Lambda}$ is identified constructively in $\mac{L}_0$ and $B$ from $(\sfc{e_n}{d_{\rm u}Y})_{n\in\Lambda}$.
\end{enumerate}
\label{thm: 1.555}
\end{cor}

\rem The same thing as mentioned in Remarks 1 and 2 in Corollary \r{thm: 1.5} holds.

\begin{prf}
We obtain (A) and (B) as in Corollary \r{thm: 1.5},
then prove (C).
We can identify $(\sfc{e_n}{a\,d_{\rm u}B})_{n\in\Lambda}$ constructively in $\mac{L}_0$ and $B$, so we identify $(\ip{e_n}{b}{})_{n\in\Lambda}$ pathwisely by $\ip{e_n}{b}{}=\sfc{e_n}{d_{\rm u}Y}-\sfc{e_n}{a\,d_{\rm u}B},\, n\in\Lambda$.
Then we identify $b^{\Lambda}\in L^0(\Omega\,;L^2([0,L]\,;\C))$
by the Fourier series expansion of its paths.
\end{prf}

\subsection{Identification of S-type It\^o processes or more general \wi functionals} 

In this subsection,
we give an additional result about identification from SFC-Os of a random function of unbounded variation.
From \tr{thm: 17.02} and \tr{thm: 21} (Theorem 4.2 in \cite{KH17b}), we have the following similar to \tr{thm: 0003} (Theorem 4.3 in \cite{KH17b}).

\begin{thm}
Let $\L_0$ be the assigned first-order language defined in Example \r{ex:1}.
Let $(e_n)_{n\in\mathbb{N}}$ and $(\vp_m)_{m\in\mab{N}}$ be CONSs of $L^2([0,1];\mab{R})$ and $a(t)$ a process defined by \eqref{eq:17} and $b\in\mathcal{L}_1^{0,2}$.
We assume the condition
\begin{equation}
\int_0^{1}K(t,s)D_th(s)\,ds\in\mac{L}_1^{0,2}
\label{eq:18}
\end{equation}
and assume that $e_n(t)$, $(\vp_m)_{m\in\mab{N}}$ and $f(t)$ satisfy any of the
conditions (C.1),(C.2) and (C.3) in Subsection \ref{sub: 1.1} for each $n\in\mab{N}$.
Then, letting
$$
d_{\dagger}Y_t:=\begin{cases}
	a(t)\,d_{\rm{u}}B_t+b(t)\,dt& \text{, if (C.3) holds for each $e_n(t)$}\\
	a(t)\,d_{\varphi}B_t+b(t)\,dt& \text{, otherwise}
	\end{cases}
$$
$a(t)$ and $b(t)$ are identified constructively in $\mac{L}_0$ and $B$ from the system $((e_n,d_{\dagger}Y))_{n\in\mab{N}}$ of SFC-$\rm{O}_{\dagger}$'s.
\label{thm: a5}
\end{thm}
\remi{1} If $\tsup_{t\in[0,1]}|K(t,\cdot)|_{L^2[0,1]}<\infty$ or $\tsup_{t\in[0,1]}|D_th|_{L^2([0,1]\t\Om)}<\infty$, \eqref{eq:18} holds.\\
\remi{2} When $K(t,s)=\ic{s\le t}$, namely,  $a(t)$ is an S-type It\^o process defined by \eqref{al:15.7}, $\tsup_{t\in[0,1]}|K(t,\cdot)|_{L^2[0,1]}<\infty$ is satisfied.
\begin{prf} 
By the definition of $(e_n,d_{\dagger}Y)$ and \tr{thm: 17.02} we have
\begin{align*}
(e_n,d_{\dagger}Y)=&
\begin{cases}
(e_n,a\,d_{\rm{u}}B)+\ip{e_n}{b}{L^2[0,1]}&\\
(e_n,a\,d_{\varphi}B)+\ip{e_n}{b}{L^2[0,1]}
\end{cases}\\
=&(e_n,a\,\de B)+\ip{e_n}{c}{L^2[0,1]},
\end{align*}
where $c(t)=\frac{1}{2}f(t)+\int_0^tD_tf(s)\,\de B_s+\int_0^1K(t,s)D_th(s)ds+D_ta(0)+b(t)\in \mac{L}_1^{0,2}$.
Then, $(e_n,d_{\dagger}Y)$ is SFC-S of $a\,\de B+c\,d\lambda$ and we can see $a(t)$ is identified constructively in $\mac{L}_0$ and $B$
from $((e_n,d_{\dagger}Y))_{n\in\mab{N}}$ by following the proof of 
Theorem 4.2(Theorem \r{thm: 21} in this note) in \cite{KH17b}.
Now, by the definition of $(e_n,d_{\dagger}Y)$ we identify $(\ip{e_n}{b}{L^2[0,1]})_{n\in\mab{N}}$ constructively in $\mac{L}_0$ and $B$, then identify $b(t)$ by the Fourier expansion of its paths.
\end{prf}


\section{Identification of random functions from \sfcsS}

In this section, we give the main result about identification of random functions from from \sfcss of a stochastic differential by $a(t)$ and $b(t)$.

\subsection{Identification of locally absolutely continuous \wifs}

Applying Theorems \r{thm:6.220}, \r{pr:4.11440} and Corollary \r{thm: 1.555} proved in Subsection 6.1, we have the following main results which give reconstructions by the law of iterated logarithm
of \wifs from \sfcsS.

\begin{thm}[Constructive identification]
Let $\Lambda$ be a cofinite subset of $\mab{N}$.
Suppose $e_n(t)$ has compact support for each $n\in\Lambda$ and
$a\in L^0([0,L]\t\Om)$ with the representative $\hat{a}(t)$ which satisfies the following:
\begin{enumerate}
\i$\hat{a}(t)$ is locally absolutely continuous in $t$ a.s.
\i$\hat{a}'(t):=\frac{d}{dt}\hat{a}(t)\in\mathcal{L}_1^{1,2},\,\hat{a}(0)\in\mathcal{L}_0^{1,2}$.
\i$\hat{a}'(t)\in L^1[0,L]$ a.s. , $\int_0^tD_t\hat{a}'(s)\,ds\in L^2[0,L]$ a.s. 
\end{enumerate}
Then, letting $\de X_t=a(t)\,\de B_t+b(t)\,dt$, $b\in L^0(\Omega\,;L^2([0,L]\,;\mab{C}))$
the following hold:
\begin{align}
\,\,\displaystyle\Biggl(\,\,|a|(t)=\lims{s\searrow t}\frac{\,\widetilde{\mathcal{P}}((e_n,\de X))_{n\in\Lambda}(s)-\widetilde{\mathcal{P}}((e_n,\de X))_{n\in\Lambda}(t)}{\sqrt{2(s-t)\log\log\frac{1}{s-t}}} \q \tx{a.s.} \,\,\Biggl)\q \tx{a.a.} \,\,t\in[0,L],
\label{al:4.11154500}
\end{align}
where $\widetilde{\mathcal{P}}$ is defined
in Lemma \ref{lem:6.3}.
Therefore, $|a|$ is $\mac{L}_0$-constructively (thus, and B-independently) identified from $((e_n,\de X))_{n\in\Lambda}$ by \eqref{al:4.11154500},
where $\L_0$ is the assigned first-order language defined in Example \r{ex:1}.
\label{thm:6.2000}
\end{thm}
\rem If $L<\infty$, then the condition 3 holds. \\

\begin{thm}[Constructive identification with $B$]
Let $\Lambda$ be a cofinite subset of $\mab{N}$.
Suppose $e_n(t)$ has compact support for each $n\in\Lambda$ and
$a\in L^0([0,L]\t\Om)$ with the representative $\hat{a}(t)$ which satisfies
the assumptions from 1 to 3 in Theorem \ref{thm:6.2000}.
Then, letting $\de X_t=a(t)\,\de B_t+b(t)\,dt$, $b\in L^0(\Omega\,;L^2([0,L]\,;\mab{C}))$
the following hold:
\begin{align}
\,\,\displaystyle\Biggl(\,\,a(t)=\lim_{k\to\infty}\Biggl(\,\lims{s\searrow t}\frac{\,\widetilde{\mathcal{P}}((e_n,\de X))_{n\in\Lambda}(s,t)+k(B_s-B_t)}{\sqrt{2(s-t)\log\log\frac{1}{s-t}}}-k\,\Biggl) \q \tx{a.s.} \,\,\Biggl)\q \tx{a.a.} \,\,t\in[0,L],
\label{al:4.11155500}
\end{align}
where $\widetilde{\mathcal{P}}$ is defined
in Lemma \ref{lem:6.3}
and $\widetilde{\mathcal{P}}((e_n,\de X))_{n\in\Lambda}(s,t)$ denotes
$\widetilde{\mathcal{P}}((e_n,\de X))_{n\in\Lambda}(s)-\widetilde{\mathcal{P}}((e_n,\de X))_{n\in\Lambda}(t)$.
Therefore, $a(t)$ is identified constructively in $\mac{L}_0$ and $B$ from $((e_n,\de X))_{n\in\Lambda}$ by \eqref{al:4.11155500},
where $\L_0$ is the assigned first-order language defined in Example \r{ex:1}.
\label{pr:4.11400}
\end{thm}

\begin{cor}
Let $\L_0$ be the assigned first-order language defined in Example \r{ex:1}
and $\Lambda$ a cofinite subset of $\mab{N}$.
Suppose $e_n(t)$ has compact support for each $n\in\Lambda$.
Let $\hat{a}(t)$ be a random function on $[0,L]$ which satisfies
the assumptions from 1 to 3 in Theorem \ref{thm:6.2000}.
Put $a=[\hat{a}]\in L^0([0,L]\t\Om)$.
Then, letting $\de X_t=a(t)\,\de B_t+b(t)\,dt$, $b\in L^0(\Omega\,;L^2([0,L]\,;\mab{C}))$,
$b^{\Lambda}=\tsum_{n\in\Lambda}\ip{e_n}{b}{}e_n\in L^0(\Omega\,;L^2([0,L]))$
the following hold:
\begin{enumerate}
\i[(A)] $|\hat{a}|$ is $\mac{L}_0$-constructively (thus, and B-independently) identified almost surely from $(\sfc{e_n}{\de X})_{n\in\Lambda}$
by \eqref{al:4.11154500} and
\begin{equation*}
|\hat{a}|(t)=
\begin{cases}
\lims{n\to\infty}n\int_t^{t+\frac{1}{n}} |a|\,d\lambda	&,t<L\\
\lims{n\to\infty}n\int_{t-\frac{1}{n}}^t |a|\,d\lambda	&,t=L<\infty.
\end{cases}
\end{equation*}
\i[(B)] $\hat{a}$ is identified constructively in $\mac{L}_0$ and $B$
almost surely from $(\sfc{e_n}{\de X})_{n\in\Lambda}$
by \eqref{al:4.11155500} and
\begin{equation*}
\hat{a}(t)=
\begin{cases}
\lims{n\to\infty}n\int_t^{t+\frac{1}{n}} a\,d\lambda	&,t<L\\
\lims{n\to\infty}n\int_{t-\frac{1}{n}}^t a\,d\lambda	&,t=L<\infty.
\end{cases}
\end{equation*}
\i[(C)] $b^{\Lambda}$ is identified constructively in $\mac{L}_0$ and $B$ from $(\sfc{e_n}{\de X})_{n\in\Lambda}$.
\end{enumerate}
\label{thm: 15}
\end{cor}

\rem Almost surely, 
$|\hat{a}|(t)$ (so is $\hat{a}(t)$) is identified at every $t\in[0,L]$.\\

We prove \tr{thm:6.2000}, \r{pr:4.11400} and Corollary \r{thm: 15}
altogether.
\begin{prf} 
For each $n\in\Lambda$, 
by Corollary \r{cor:13.7} $\overline{e_n}a$ is u-integrable for $L^2([0,L]\,;\mab{C})$ and the \og integral is given by
$$
\int_0^L\overline{e_n(t)}a(t)\,d_{\rm u}B_t=\int_0^L\overline{e_n(t)}a(t)\,\de B_t+\int_0^{L}\overline{e_n(t)}\Bigl(\,\int_0^{t}D_ta'(s)\,ds+D_ta(0)\,\Bigl)\,dt.
$$
Thus, 
$\sfc{e_n}{\de X}=\int_0^L \overline{e_n(t)}a(t)\,d_{\rm u}B_t+\int_0^L \overline{e_n(t)}c(t)\,dt$,
where $c(t)=b(t)-\int_0^{t}D_ta'(s)\,ds-D_ta(0)$.
Because $c\in L^2([0,L]\,;\mab{C})$ almost surely from the assumptions,
$\sfc{e_n}{\de X}=\sfc{e_n}{a\,d_{\rm u}B}+\ip{e_n}{c}{}$ is \sfco of $a\,d_{\rm u}B
+c\,d\lambda$. In addition, $a(t)$ is a \nc finite variation process from the assumptions (1) and $a'\in L^1[0,L]$ almost surely.
Thus by Theorem \r{thm:6.220}, \r{pr:4.11440} and (A),(B) in Corollary \r{thm: 1.555}
we obtain \tr{thm:6.2000}, \r{pr:4.11400} and (A),(B) in Corollary \r{thm: 15}, respectively. 
Now that $a(t)$ is identified, we can identify $\ip{e_n}{b}{}=\sfc{e_n}{\de X}-(e_n,a\,\de B)$ constructively in 
$\mac{L}_0$ and $B$ from $(\sfc{e_n}{\de X})_{n\in\Lambda}$
for each $n\in\Lambda$. Then, we identify $b^{\Lambda}\in L^0(\Omega\,;L^2([0,L]))$ 
by the Fourier series expansion of its paths.
\end{prf}





\section{Appendix A: Notion of constructiveness}
\label{sub:3.2}

In Appendixes A and B,
we follow the notation and terminology in Subsection \r{sub:011}.

In this section, we introduce the notion of constructiveness in an assigned first-order language.
Note that the notion is metamathematical.
For the fundamental notion of first-order logic, refer to \cite{En} or \cite{Sh}.
Let $(e_n)_{n\in\mab{N}}$ be the CONS fixed at the beginning of Section \ref{sec:4}.
Let $\L$ be a language of first-order logic
and
$V=\{x\,|\,x=x\,\}$ denotes the ZFC universe, i.e.
the domain of discourse of ZFC set theory.
Refer to \cite{Kn} for ZFC set theory.
The logical terms and logical formulas of $\mac{L}$ are defined by recursion on natural numbers.
Assume that there is a one-to-one correspondence between $\L$ and some set $V_0$ in $V$ as follows${}^{(\ast 1)}$:

\begin{enumerate}
\i Each constant symbol $c$ in $\L$ is assigned to an element $c^V$ in $V_0$.
\i Each $n$-place function symbol $f$ in $\L$ is assigned to an $n$-ary function $f^V$ in $V_0$. 
\i Each $n$-place relation symbol $p$ in $\L$ is assigned to an $n$-ary relation $p^V$ in $V_0$.
\end{enumerate}

\footnotetext{\normalsize $^{(\ast 1)}$
This statement seems to be the statement of the definition of a structure for the language $\mac{L}$ (see Chapter 2.5 in \cite{Sh} or Chapter 2.2 in \cite{En}).
But $V$ is the structure for the language $\{=,\in\}$  with predicate symbols '$=$' and '$\in$'  assigned to the equality relation and membership relation, respectively.
We never intend to define a model for $\mac{L}$ of some axiomatic system.
}

\begin{dfn}
We call
a first-order language with an assignment as described above  
an assigned first-order language.
\end{dfn}
In what follows, we introduce two kinds of definitions of constructiveness via evaluations of a logical term and formula of the assigned language $\mac{L}$.\m

\ni
$\bullet\,$ \textbf{Semantic definition}

\begin{dfn}[Semantic evaluation]
We define the (semantic) evaluation$^{(\ast 2)}$ ${\rm eval}_{V}(s)$ in $V$ of a closed logical term $s$ of $\mac{L}$ as follows:
\begin{enumerate}
\i If $s$ is a constant symbol, then ${\rm eval}_{V}(s)=({\rm eval}_{V}(s)_1,{\rm eval}_{V}(s)_2):=(s^V,1)$.
\i The case that $s$ is represented as $fs_1...s_n$ with terms $s_1,...,s_n$ and an $n$-place function symbol $f$: 
If $({\rm eval}_{V}(s_1)_2,...,{\rm eval}_{V}(s_n)_2)=(1,...,1)$
and 
$y=$
$f^V({\rm eval}_{V}(s_1)_1,...,$
${\rm eval}_{V}(s_n)_1)$ for some $y$,
define ${\rm eval}_{V}(s):=(y,1)$, 
otherwise
define ${\rm eval}_{V}(s):=(0,0)$.
\end{enumerate}
\end{dfn}

\footnotetext{\normalsize $^{(\ast 2)}$
This definition is same as that of interpretation or individual (in Chapter 2.5 in \cite{Sh}) in $V$ of a term of $\mac{L}$, we considering $V$ as the structure for $\mac{L}$.
The only difference is that the case $\,f^V({\rm eval}_{V}(s_1)_1,...,$
${\rm eval}_{V}(s_n)_1) \tx{ isn't defined}$ could occur in the recursive case.
Then, we add the binary variable to the evaluation amount which
returns 0 if the evaluation is ill-defined.
}

Similarly, we can also define the evaluations (evaluation and truth value) ${\rm eval}_{V}^{\sigma}$ in $V$ of an unclosed logical term, i.e. a logical term with at least one variable and of a logical formula, for each variable assignment $\sigma$, respectively (see p. 83 in \cite{En}).
\begin{dfn}[Semantically constructive map]
We say a map $\mac{T}$ is semantically constructive (or direct-representable) in $\mac{L}$ (or $\mac{L}$-semantically constructive) if there exists a logical term $s(x)$ of $\mac{L}$ with one  variable $x$ such that 
'${\rm eval}_{V}^{\sigma}(s(x))=(\mac{T}({\rm eval}_{V}^{\sigma}(x)_1),1)$ for every variable assignment $\sigma:\{x\}\to{\rm dom}\mac{T}$'.
In this case, we also say $\mac{T}$ is $\mac{L}$-semantically constructive with the term $s(x)$
and say $\mac{T}(X)$ is $\mac{L}$-semantically constructive, for each $X\in {\rm dom}\mac{T}$.
\end{dfn}\m

\ni
Here, the statement
'${\rm eval}_{V}^{\sigma}(s(x))=(\mac{T}({\rm eval}_{V}^{\sigma}(x)_1),1)$ for every variable assignment $\sigma:\{x\}\to{\rm dom}\mac{T}$'
in this definition means
'the evaluation of the corresponding formula of $\{=,\in\}$ is true in $V$',
and the evaluations of a term and formula do not mean the syntax derived values by natural deduction but the semantic assigned values in $V$.
But, one can prove in ZFC finite syntax evaluations equal the corresponding semantic evaluations since we can choose the model $R(\gamma)$ for a sufficiently large ordinal $\gamma$ as in Exercise II.18.15 in \cite{Kn} of the finite subclass of ZFC
whose axioms were used in the proof, because of the reflection theorem (see Exercise II.18.15 in \cite{Kn}).

This definition is natural as a purely mathematical notion.
On the other, there may be $\mac{L}$-constructive maps with represented formulas of $\mac{L}$
which is unprovable from ZFC.
Then, in what follows we also define an $\mac{L}$-constructive map in more strong sense,
that is to say, a constructive map
with a provable (demonstrable) formula of $\mac{L}$.\m


\ni
$\bullet\,$ \textbf{Syntactic definition}

\begin{dfn}[Syntactic evaluation]
We define the syntax evaluation ${\rm syn}_{V}(s)$ in $V$ of a closed logical term $s$  of $\mac{L}$ as follows:
\begin{enumerate}
\i If $s$ is a constant symbol, then ${\rm syn}_{V}(s)=({\rm syn}_{V}(s)_1,{\rm syn}_{V}(s)_2):=(s^V,1)$.
\i The case that $s$ is represented as $fs_1...s_n$ with terms $s_1,...,s_n$ and an $n$-place function symbol $f$: 
If $({\rm syn}_{V}(s_1)_2,...,{\rm syn}_{V}(s_n)_2)=(1,...,1)$ and 
''$y=$
$f^V({\rm syn}_{V}(s_1)_1,...,$
${\rm syn}_{V}(s_n)_1)$'' 
is provable for some $y$,
define ${\rm syn}_{V}(s):=(y,1)$,
otherwise
define ${\rm syn}_{V}(s):=(0,0)$.
\end{enumerate}
\end{dfn}\m

Similarly, we can also define the syntactic evaluations (evaluation and truth value) ${\rm syn}_{V}^{\sigma}$ in $V$ of an unclosed logical term, i.e. a logical term with at least one variable and of a logical formula,
for each variable assignment $\sigma$, respectively.

\begin{dfn}[Syntactically constructive map]
We say a map $\mac{T}$ is (syntactically) constructive (or direct-representable) in $\mac{L}$ (or $\mac{L}$-(syntactically) constructive) if there exists a logical term $s(x)$ of $\mac{L}$ with one  variable $x$ such that
''${\rm syn}_{V}^{\sigma}(s(x))=(\mac{T}({\rm syn}_{V}^{\sigma}(x)_1),1)$ for any variable assignment $\sigma:\{x\}\to{\rm dom}\mac{T}$'' is provable.
In this case, we also say $\mac{T}$ is $\mac{L}$-(syntactically) constructive with the term $s(x)$
and say $\mac{T}(X)$ is $\mac{L}$-(syntactically) constructive, for each $X\in {\rm dom}\mac{T}$.
\end{dfn}

\rem If a map $\mac{T}$ is syntactically constructive, $\mac{T}$ is constructive in both senses.
Then, $\mac{T}$ is simply said to be constructive, usually.\m

By $\mac{L}^W$ we denote the assigned first-order language whose symbols are the symbols of $\mac{L}$ and the function symbol assigned to a Brownian motion $(W_t)_{t\in[0,\infty)}$ on $(\Om,\F,P)$,
and an $\mac{L}^W$-semantically (resp. syntactically) constructive map is also said to be semantically 
(resp. syntactically) constructive in $\mac{L}$ and (or with) $W$.\bigskip

We notice that both kinds of notions of constructiveness for a map depend (only) on
what assigned language is selected.
Now, we are to give examples related in this note.
We call the following spaces the first-kind spaces regarding $([0,L],\Om)$:
\begin{enumerate}
\i $\C^d,\, L^0([0,L]^d),\, L^0(\Om)$ and $L^0([0,L]^d\t\Om)$ for each $d\in\N$,
\end{enumerate}
and call the following spaces the second-kind spaces regarding $([0,L],\Om)$:
\begin{enumerate}
\i[2.] $X^{\Z^{j}\t[0,L]^{k}\t\C^{l}}=\{f\,|\,f:\Z^{j}\t[0,L]^{k}\t\C^{l}\to X\}$
for each first-kind space $X$ and $(j,k,l)\in\N_0^3\diagdown{\{(0,0,0)\}}$.
\end{enumerate}

\begin{exm}

\ni
We define the assigned language $\L_0$ of first-order logic by the following list of symbols and assignments:
Here, each mathematical object to which each symbol is assigned is noted in parentheses.
\m

\noindent
\textbf{Constant symbols:}\m

\noindent
$\centerdot$ $0,1,2,\ldots$ ($0,1,2,\ldots\in\N_0$), $L$ (the constant $L$)\m 

\noindent
$\centerdot$ $(e_n)_{n\in\N}$ (the CONS $(e_n)_{n\in\N}\in (L^2[0,L])^{\N}$ of $L^2[0,L]$)\m


\noindent
\textbf{Function symbols:}\m

\noindent
(Functions on $\C^1,\C^2,\ldots$ 
for each $d=1,2,\ldots$, respectively)$^{(\ast 3)}$


\footnotetext{\normalsize $^{(\ast 3)}$
We intend that $\mathrm{id}_{[0,L]^1}$, $\mathrm{id}_{[0,L]^2},\ldots$ are symbols
but $\mathrm{id}_{[0,L]^c}$ with a character $c$ is not a symbol,
for instance.
We used $d$ as a meta symbol, here.
}

\noindent
$\centerdot$ $\mathrm{id}_{\C^d}$ (identity),
$\mathrm{id}_{[0,L]^d}$ (identity on $[0,L]^d$),
$\mathrm{id}_{\Z^d}$ (identity on $\Z^d$)

\noindent
$\centerdot$  $\ic{}^{\C^d}$ (indicator),
$\ic{}^{[0,L]^d}$ (indicator on $[0,L]^d$),
$\ic{}^{\Z^d}$ (indicator on $\Z^d$)

\noindent
$\centerdot$ $P^Y_1, P^Y_2,\ldots,P^Y_d$ (canonical projections on $Y=\C^d$, $[0,L]^d$, $\Z^d$),
${\rm Re}, {\rm Im}$ (real, imaginary parts on $\C$, respectively)

\noindent
$\centerdot$
$\inc_{\R^d, \C^d}$ (inclusion from $\R^d$ to $\C^d$),
$\inc_{\N^d, \Z^d}$ (inclusion from $\N^d$ to $\Z^d$),
$\inc_{\C^d, \C^k}$ (inclusion from $\C^d$ to $\C^k$ for each $k=d+1,d+2,\ldots$),
$\inc_{[0,L]^d, [0,L]^k}$ (inclusion from $[0,L]^d$ to $[0,L]^k$ for each $k=d+1,d+2,\ldots$),
$\inc_{\Z^d, \Z^k}$ (inclusion from $\Z^d$ to $\Z^k$ for each $k=d+1,d+2,\ldots$)



\noindent
$\centerdot$ {\rm lim} (limit in Euclid norm on $\R$)


\noindent
$\centerdot$ {\rm sup}, {\rm inf} (supremum, infimum on $\R$, respectively)


\m\noindent
(Operators on  $L^2[0,L]$)

\noindent
$\centerdot$ $J_1,J_2,\ldots$ (orthogonal projections onto $e_1,e_2,\ldots$, respectively)

\m\noindent
(Operations on  $L^0([0,L]^1), L^0([0,L]^2),\ldots$ for each $d=1,2,\ldots$, respectively)

\noindent
$\centerdot$ $\int d\lambda^{\ot d}$ (integral with respect to $\lambda^{\ot d}$ if it is defined)

\noindent
$\centerdot$ ${\rm lim}_{*}$
(limits in measure, in $L^p$ (for each $p\in[1,\infty]$) and with measure 1,
with respect to $\lambda^{\ot d}$)

\m\noindent
(Operations on  $L^0(\Om)$)

\noindent
$\centerdot$ $E$ (expectation with respect to $P$),
$\int dP\ot\lambda^{\ot d},\int d\lambda\ot P\ot\lambda^{\ot d-1},\ldots,\int d\lambda^{\ot d}\ot P$
(integrals with respect to $P\ot\lambda^{\ot d},\lambda\ot P\ot\lambda^{\ot d-1},\ldots,\lambda^{\ot d}\ot P$)

\noindent
$\centerdot$ ${\rm lim}_{*}$
(limits in probability, in $L^p(\Om)$ and with probability 1, and
in measure, in $L^p$ and with measure 1 with respect to $P\ot\lambda^{\ot d},\lambda\ot P\ot\lambda^{\ot d-1},\ldots,\lambda^{\ot d}\ot P$)

\m\noindent
(Operations$^{(\ast 4)}$ on each first or second-kind space $X$)

\footnotetext{\normalsize $^{(\ast 4)}$
The symbols which are assigned to these operations depend on $X$.
But one can extend the domain of each of these operations independently of $X$.  
Thus, we can omit the index $X$ from each symbol.
}

\noindent
$\centerdot$ 
The corresponding symbols to the above symbols
(the corresponding functions on $X^{\C^d}$
to the above functions on $\C^d$
)$^{(\ast 5)}$

\footnotetext{\normalsize $^{(\ast 5)}$
For instance, the corresponding function to the inclusion from $\R^d$ to $\C^d$
is the inclusion from $X^{\R^d}$ to $X^{\C^d}$.
}

\noindent
$\centerdot$
The corresponding symbols to the above symbols
(the corresponding functions on the second-kind space $X^{\Z^{j}\t[0,L]^{k}\t\C^{l}}$ to the above functions on the first-kind space $X=\C^d$, $L^2[0,L]$, $L^0([0,L]^d)$, $L^0(\Om)$ )$^{(\ast 6)}$

\footnotetext{\normalsize $^{(\ast 6)}$
For instance, the corresponding function (operator) to the integral with respect to $\lambda^{\ot d}$
is the integral with respect to $\lambda^{\ot d}$ with parameter $x\in\Z^{j}\t[0,L]^{k}\t\C^{l}$.
}

\noindent
$\centerdot$ $+,-,\times,/$ (canonical sum, subtraction, multiplication and division (on the domains), respectively),
$\cdot\,$(scalar multiplication)

\noindent
$\centerdot$
$\tsum$ ($X^{\Z}\to X^{\Z}$ which associates
$(x_n)_{n\in\Z}$ with $(\,\sum_{-n\le i\le n,\,\dot{\vp}(n)}x_i\,)_{n\in\Z}$ (partial sum))

\noindent
$\centerdot$
$\tprod$
($X^{\Z}\to X^{\Z}$ which associates
$(x_n)_{n\in\Z}$ with $(\,\prod_{-n\le i\le n,\,\dot{\vp}(n)}x_i\,)_{n\in\Z}$ (partial product))
(we noted below for $\dot{\vp}(n)$)

\noindent
$\centerdot$
$s(\dot{x})$
($X\to Y$ which associates $a\in X$ with ${\rm eval}_{V}^{\sigma_a}(s(x))_1\in Y$,
where $\sigma_a:\{x\}\to\{a\}$,
for each term $s(x)$ of $\L_0$ with a variable $x\in X$)

\noindent
$\centerdot$
$P^{\Z},\,P^{[0,L]},P^{\C}$ (canonical projections with parameters $n\in\Z,t\in[0,L],x\in\C$ from $X^{\Z},X^{[0,L]},X^{\C}$ to $X$, respectively)

\noindent
$\centerdot$
$J^{\Om}$ (canonical projection with parameter $\om\in\Om$ from $X$ to the space removed $\Om$ from $X$
$^{(\ast 7)}$, where $X$ is $L^0(\Om)$ or $L^0([0,L]^d\t\Om)$ or the second-kind space for $L^0(\Om)$ or $L^0([0,L]^d\t\Om)$.

\footnotetext{\normalsize $^{(\ast 7)}$
For instance if $X=L^0(\Om)^{[0,L]}$,
the space removed $\Om$ from $X$ is $\C^{[0,L]}$.
}

\noindent
$\centerdot$
$J^{[0,L]}$ (canonical projection with parameter $t\in[0,L]$ from $X$ to the space removed $[0,L]$ from $X$, 
where $X$ is $L^0([0,L]^d)$ or $L^0([0,L]^d\t\Om)$ or the second-kind space for $L^0([0,L]^d)$ or $L^0([0,L]^d\t\Om)$.

\noindent
$\centerdot$
$\inc^{\Z},\,\inc^{[0,L]},\,\inc^{\C}$ (inclusions from $X$ to $X^{\Z},\,X^{[0,L]},\,X^{\C}$, respectively)


\noindent
$\centerdot$
$\inc_{\C^d, L^0(\Om)^d}$ (inclusion from $\C^d$ to $L^0(\Om)^d$),
$\inc_{\C^{[0,L]^d}, L^0([0,L]^d)}$ (inclusion from the measurable subset of $\C^{[0,L]^d}$ to $L^0([0,L]^d)$) 
and the other symbols for the same-type canonical inclusions defined\m

\m
\noindent
\textbf{Predicate symbols:}
$=$ (equality in $V$), $<$ (inequality in $\R$)
\m

Here indicator, limit, supremum, infimum, partial sum and partial product are determined for each formula of $\L_0$ without quantifier  with respect to the parameter set, $\Z^d$, $[0,L]^d$ or $\C^d$.


\ni
Then, the following hold:
\begin{enumerate}
\i[(a)] $e,\pi,\sqrt{-1}$ are constructive in $\L_0$.
\i[(b)] The elementary functions: $\sin, \cos, \tan,\arcsin, \arccos, \arctan,{\rm pow}, {\rm log}$ (on the domains) are constructive in $\L_0$.
\i[(c)] $f(x)=\ic{\Z}(x),\,x\in\R$ is constructive in $\L_0$.
\i[(d)] $g(x)=\ic{\Q}(x),\,x\in\R$ (Dirichlet function) is constructive in $\L_0$.
\i[(e)] $\mathscr{T}_{\tau, T}((\hat{a}_n)_{n\in\mab{Z}})=\hat{a}_0+\tsum_{n\neq 0}\frac{1}{2\pi\sqrt{-1}n}\hat{a}_n T_n$, where $T_n$ and $\hat{a}_n$ is defined in Theorem \r{thm: 0001}, is constructive in $\L_0$.
\i[(f)] $\mathscr{F}_{\tau, T}$ in Theorem \r{thm: 0001} is constructive in $\L_0$ with $B$.
\end{enumerate}
\label{ex:1}
\end{exm}

\begin{prf}
(a) and (b) are obvious.
(c) and (d) hold since $f(x)=\tlim_{n\to\infty}\cos^{2n}(\pi x)$ and
$g(x)=\tlim_{m\to\infty}f(m!\,x)$, as are well known.
(e) holds because the exponential functions $T_n,\,n\in\Z$ are $\L_0$-constructive
and $\neg\,\, n= 0$ is the formula of $\L_0$
with respect to $\Z$.
(f) holds because the SFCs $\hat{a}_n,\,n\in\mab{Z}$ are constructive in $\L_0$ and $B$,
for the Ogawa integrals with respect to the Haar system are constructive in $\L_0$ and $B$.
\end{prf}






\section{Appendix B}

The following Propositions \r{ap:21} and \r{ap:24} are concerned with measurability and continuity of a stochastic process whose paths are of bounded variation. 

\begin{prp} 
Let $a:[0,L]\t\Om\to\mab{R}$ be a weak stochastic process, namely, $a_t$ is a random variable for every $t\in[0,L]$.
Suppose all paths of $a(t)$ are left-continuous and of bounded variation.
Then, both $a_{+}(t)$ and $a_{-}(t)$ are weak stochastic processes,
and their all paths are left-continuous.
\label{ap:21}
\end{prp}
\rem In particular, both $a_{+}(t)$ and $a_{-}(t)$ become measurable stochastic process since any weak stochastic process whose all paths are left-continuous is measurable (see Proposition 1.13 in \cite{KS} for example).
\begin{prf} 
First, the left-continuity of $a_{+}(t)$ and $a_{-}(t)$ follow the left-continuity of $a(t)$.
Next, we show that $a_{-}(t)$ is weak stochastic process (we can show so is $a_{+}(t)$). Fix any $t\in(0,L]$.
Put $V(t_0,...,t_n)=\tsum_{j=1}^{n}(a(t_j)-a(t_{j-1}))^{-}$ for $n\in\mab{N}$ and $t_0,\ldots,t_n$ such that $0=t_0<t_1<\cdots<t_n=t$.
All of them are random variables since $a(t)$ is the weak stochastic process.
While $a_{-}(t)$ is represented as
\begin{align}
a_{-}(t)=\sup_{\substack{n\in\mab{N}	\\ 0=t_0<t_1<\cdots<t_n=t}}V(t_0,...,t_n),
\label{al:21.5}
\end{align}
setting $A=\{\,(t_0,...,t_n)\,|\,n\in\mab{N},\, 0=t_0<t_1<\cdots<t_n=t,\, t_1,...,t_{n-1}\in\mab{Q}\,\}$, \eqref{al:21.5} is rephrased as
$$
a_{-}(t)=\sup_{(t_0,...,t_n)\in A}V(t_0,...,t_n)
$$
by the denseness of $\mab{Q}$ and the left-continuity of $a(t)$.
Besides, the cardinality of $A$ equals that of $\d\cu_{n=1}^{\infty}\mab{Q}^{n-1}$ ($\mab{Q}^{0}$ denotes some singleton), which is a countable set. Therefore, $a_{-}(t)$ is a random variable.
\end{prf}

\begin{prp}  
Let $a(t)$ be a \nc finite variation process on $[0,L]$.
Then, there exists a version $\widetilde{a}(t)$ of $a(t)$ in $L^0([0,L]\t\Om)$ such that
\begin{enumerate}
\i All paths of $\,\widetilde{a}(t)$ are left-continuous and of bounded variation,
\i $\widetilde{a}_{+}(t)$ and $\widetilde{a}_{-}(t)$ are $\L([0,L])\ot\F$-measurable and all their paths are left-continuous and monotonically increase.
\end{enumerate}
\label{ap:24}
\end{prp} 
\begin{prf} 
By the assumption on $a(t)$, there exists a probability-1 set $\Om^V\in\F$ such that $a(t)$ is of bounded variation on $\Om^V$.
Set $a^{(0)}=a\ic{\Om^V}$, then $a^{(0)}(t)$ is a version of $a(t)$ in $L^0([0,L]\t\Om)$ and all paths of $a^{(0)}(t)$ are of bounded variation.
Furthermore, set $\widetilde{a}(t,\om)=\tlim_{s\nearrow t}a^{(0)}(s,\om)$, which is $\L([0,L])\ot\F$-measurable.
Since discontinuous points of a function of bounded variation are at most countable, 
$\widetilde{a}(t)$ is a version of $a^{(0)}(t)$ in $L^0([0,L]\t\Om)$ and all paths of  $\,\widetilde{a}(t)$ are left-continuous and of bounded variation.
Therefore by Proposition \r{ap:21},  $\widetilde{a}_{+}(t)$ and  $\widetilde{a}_{-}(t)$ are measurable stochastic processes whose all paths are left-continuous and monotonically increase.
\end{prf}



In what follows, using Proposition \r{ap:21}, we show that 
$$
S=\{(t,\om)\in[0,L]\t\Om\,|\,a_{+}(\cdot,\om)\tx{ and }a_{-}(\cdot,\om)\tx{ are differentiable in $t$}\}
$$
is a measurable set of $[0,L]\t\Om$ for a (\nC) stochastic process $a:[0,L]\t\Om\to\mab{R}$ whose all paths are left-continuous and of bounded variation.

\begin{prp} 
Let $a:[0,L]\t\Om\to\mab{R}$ be $\L([0,L])\ot\F$-measurable.
Assume all paths of $a(t)$ are bounded and they monotonically increase.
Then, all Dini derivatives of $a(t)$
$$
D_{\pm}a(t,\om)=\limi{h\to\pm0}\frac{a(t+h,\om)-a(t,\om)}{h}
,\q
D^{\pm}a(t,\om)=\lims{h\to\pm0}\frac{a(t+h,\om)-a(t,\om)}{h}
$$
are $\L([0,L])\ot\F$-measurable.
\label{lem:20}
\end{prp}
\begin{prf} 
Put $g_h(t,\om)=\tinf_{0<|\delta|<|h|}\frac{a(t+\delta,\om)-a(t,\om)}{\delta}$ for each $h\in(-1,1)$.
We show the measurability of $D_{+}a(t,\om)$ (the others are shown by similar arguments).
Since
$
D_{+}a(t,\om)
=\tlim_{n\to\infty}g_{\frac{1}{n}}(t,\om),
$
we have only to show for any $h>0\,$ $g_h$ is $\L([0,L])\ot\F$-measurable, in other words,
$g_h^{-1}(\{-\infty\}\cup(-\infty,c))\in\L([0,L])\ot\F$ for any real number $c$.
At first, we have
\begin{align*}
g_h^{-1}(\{-\infty\}\cup(-\infty,c))
&=\{\,(t,\om)\in[0,L]\t\Om\,|\,\e \delta\in(0,h)\,\,\frac{a(t+\delta,\om)-a(t,\om)}{\delta}<c\}.
\end{align*}
Now,  we show that for any $(t,\om)\in[0,L]\t\Om$
\begin{align}
\e \delta\in(0,h)\,\,\frac{a(t+\delta,\om)-a(t,\om)}{\delta}<c\q \Rightarrow\q \e r\in(0,h)\cap\mab{Q}\,\,\frac{a(t+r,\om)-a(t,\om)}{r}<c.
\label{al:20.55}
\end{align}
Assume the antecedent and fix one of the $\delta\in(0,h)$. for the denseness of $\mab{Q}$ there exists $r\in(0,\delta)\cap\mab{Q}$ such that $\frac{a(t+\delta,\om)-a(t,\om)}{r}<c$. Here, because $a(\cdot,\om)$ monotonically increase we have
$\frac{a(t+r,\om)-a(t,\om)}{r}\le\frac{a(t+\delta,\om)-a(t,\om)}{r}<c$. 
Thus,  the consequent holds. 
Also, the converse of \eqref{al:20.55} is trivial, so we have the following:
\begin{align*}
g_h^{-1}(\{-\infty\}\cup(-\infty,c))
&=\{\,(t,\om)\in[0,L]\t\Om\,|\,\e r\in(0,h)\cap\mab{Q}\,\,\frac{a(t+r,\om)-a(t,\om)}{r}<c\}\\
&=\cu_{r\in(0,h)\cap\mab{Q}}\Bigl(\,\frac{a(\cdot+r)-a(\cdot)}{r}\,\Bigl)^{-1}((-\infty,c)).
\end{align*}
$\frac{a(\cdot+r)-a(\cdot)}{r}$ is $\L([0,L])\ot\F$-measurable function for $r$ since 
$a(t)$ is $\L([0,L])\ot\F$-measurable, then $g_h^{-1}(\{-\infty\}\cup(-\infty,c))\in\L([0,L])\ot\F$ is proved.
\end{prf}

\begin{prp} 
Let $a:[0,L]\t\Om\to\mab{R}$ be a  (\nC) stochastic process whose all paths are left-continuous and of bounded variation.
Then,
$S=\{(t,\om)\in[0,L]\t\Om\,|\,$
$a_{+}(\cdot,\om)\tx{ and }a_{-}(\cdot,\om)\tx{ are differentiable}$
$\tx{ in $t$}\}$
is a measurable set of $[0,L]\t\Om$. In particular, the following are equivalent.
\begin{enumerate}
\i[$(i)\,$] $a(\cdot,\om)$ and $a_{\rm tv}(\cdot,\om) \tx{ are differentiable in } t \,\tx{ a.a.} \,\, (t,\om)\in[0,L]\t\Om$.
\i[$(ii)\,$] \mbox{$\tx{( } a(\cdot,\om)$ and $a_{\rm tv}(\cdot,\om) \tx{ are differentiable in } t \,\tx{ a.a.} \,\, t\in[0,L] \tx{ ) \, a.a.} \,\, \om\in\Om$.}
\i[$(iii)$] $\tx{( } a(\cdot,\om)$ and $a_{\rm tv}(\cdot,\om) \tx{ are differentiable in } t \,\tx{ a.a.} \,\, \om\in\Om \tx{ ) \, a.a. } \,\, t\in[0,L]$.
\end{enumerate}
\label{ap:24.5}
\end{prp} 
\begin{prf} 
For convenience, we denote $a_{+}$ and $a_{-}$ by $a_1$ and $a_2$, respectively. 
For $i\in\{1,2\}$ set
$$
D_{\pm}a_i(t,\om)=\limi{h\to\pm0}\frac{a_i(t+h,\om)-a_i(t,\om)}{h}
,\q
D^{\pm}a_i(t,\om)=\lims{h\to\pm0}\frac{a_i(t+h,\om)-a_i(t,\om)}{h}.
$$
Then, by Proposition \r{ap:21} each $a_i$ turns out to be $\L([0,L])\ot\F$-measurable, thus by Proposition \r{lem:20},
$S=\d\ca_{i=1,2}\{\,(t,\om)\in[0,L]\t\Om\,|\,	D^{+}a_i(t,\om)\le D_{-}a_i(t,\om)\tx{ and }
D^{-}a_i(t,\om)\le D_{+}a_i(t,\om)\tx{ and }$
$|D^{+}a_i(t,\om)|<\infty\,\}$
belongs to $\L([0,L])\ot\F$.
The latter statement holds 
for
$S=\{(t,\om)\in [0,L]\t\Om\,|\,a(\cdot,\om)\tx{ and }$
$a_{\rm tv}(\cdot,\om)\tx{ are differentiable in }t \}.$
\end{prf}

\bigskip
\par\noindent
{\bf Acknowledgments}
The author would like to express his grate gratitude to Professor Tetsuya Kazumi for supervising
to draw this paper in English and to Professor Shigeyoshi  \og and Professor Hideaki \ue for paying kind attention and giving me some advice and comments to this study on several occasions.
I would also like to thank Professor Masaru Kada for supervising me in my undergraduate years
and proposing literatures in the field of mathematical logic related to this paper.


\end{document}